\documentclass{article}
\usepackage{fullpage,amsmath,amssymb}
\usepackage{natbib}

\usepackage[ruled,vlined,linesnumbered]{algorithm2e} 
\usepackage{graphicx,subfig}

\SetArgSty{rm}
\SetKwInOut{Output}{Output} 
\SetKw{KwAnd}{and}
 
\usepackage{pgfplots}
\usepackage{listings}

\usetikzlibrary{plotmarks}
\newlength\figureheight 
\newlength\figurewidth   
\pgfplotsset{compat=newest}    

\usepackage{supertabular}

\renewcommand\cite[1]{\citet{#1}}

\newif\iffinal % introduce a switch for draft vs. final document
\finaltrue % use this to compile the final document

\graphicspath{{Figures/}}

\iffinal
  
\else                           
  \pgfrealjobname{paper}
  
\fi   
\newcommand{\p}{\partial}
\newcommand{\one}{{\mathbf 1}}    
    
\newcommand{\nw}{{n_\text{w}}}

\newcommand{\define}{\overset{\text{def}}{=}}
                                   
\usepackage{xcolor}  
\begin{document} 

\Large

\begin{center}
	{\bf Computing reaction rates in bio-molecular systems using discrete macro-states} \\[1em]
	Eric Darve$^{a,b}$ and Ernest Ryu$^a$ \\[1em]
	$^a$ Institute for Computational and Mathematical Engineering \\
	$^b$ Mechanical Engineering Department \\
	Stanford University \\[1em]
	Eric Darve, corresponding author: darve@stanford.edu \\[1em]
	Eric Darve, 496 Lomita Mall, Stanford CA 94305 \\
	Ernest Ryu, 496 Lomita Mall, Stanford CA 94305 \\[1em]
	\today
\end{center}

\normalsize

\clearpage

\tableofcontents

\clearpage

\section{Introduction}

Computing reaction rates in biomolecular systems is a common goal of molecular dynamics simulations. The reactions considered often involve conformational changes in the molecule, either changes in the structure of a protein or the relative position of two molecules, for example when modeling the binding of a protein and ligand. Here we will consider the general problem of computing the rate of transfer from a subset $A$ of the conformational space $\Omega$ to a subset $B \subset \Omega$. It is assumed that $A$ and $B$ are associated with minimum energy basins and are long-lived states. See Table~\ref{tab1} for the notations used in this paper.

Rates can be obtained using many different methods. In this paper we will review some of the most popular approaches. We organize the different approaches roughly in chronological order and under four main categories: reactive flux, transition path sampling, conformation dynamics. The fourth class of methods, to which we do not give any specific name, in some sense attempts to combine features from transition path sampling and conformation dynamics. They include non-equilibrium umbrella sampling (\cite{Warmflash:2007dz, Dickson:2009gt}), and weighted ensemble dynamics (\cite{Huber:1996dn}).

\medskip

{\bf Reactive flux.} We start with methods that were derived around 1930s (\cite{Marcelin:1915uv, Eyring:1931wh, Eyring:1935dt, Horiuti:1938ti, Wigner:1938uw}), were revisited later by, for example, \cite{Chandler:1978hh}, and are based on the concept of reactive flux. In these methods, the rate is derived from the free energy, and it is assumed that it is controlled by the flux at a saddle point at the top of the energy barrier separating $A$ and $B$. The advantages of this approach is that it involves quantities which are relatively easy to calculate. However it makes relatively strong assumptions about the system, and in practice assumes that a lot of information is already available regarding the transition mechanism and important pathways between $A$ and $B$. We relate some of the results to Kramers' method, which applies to systems modeled using Langevin dynamics and overdamped dynamics (\cite{Gardiner:1997tb,Hanggi:1990en}).

\medskip

{\bf Transition path sampling.} Many of the ideas developed in that context were used to develop another class of methods based on sampling transition pathways between $A$ and $B$ (\cite{Bolhuis:2002ew,1998JChPh.108.1964D,Dellago:2002uf}). From the ensemble of pathways, rates and other properties can be obtained. The advantage of some of these approaches is that they do not require determining the saddle point separating $A$ and $B$, and they apply to more general situations, for example when multiple pathways contribute to the rate. See~\cite{E:tb} for a discussion of transition-path theory, which proposes a mathematical framework to study transition pathways and the associated probability current. To address shortcomings of some of these approaches, other methods were pursued along similar lines, including transition interface sampling (\cite{vanErp:2003ik}), and forward flux sampling (\cite{Allen:2005dn, Allen:2006ch, Allen:2006cp}). We will present in the same category the milestoning technique, \cite{Faradjian:2004cr}, which although different in spirit, shares some similarities with transition interface sampling. This approach requires that the system ``loses'' memory when moving from a milestone to the next, for example by ensuring that the milestones are sufficiently separated from one another. Recent advances will be presented based on the work of \cite{VandenEijnden:2008bn}, who introduced the concept of optimal milestones using the committor function. Using these milestones, the rate can be obtained exactly (in the absence of statistical errors), even when the milestones are close to one another. \cite{Majek:2010kq} attempted to define the milestones in a way that would be computationally more general and advantageous compared to the original formulation.

\medskip

{\bf Conformation dynamics.} This is a large class of methods that can be traced back to Deuflhard and Sch\"utte (\cite{Deuflhard1998,Schutte2003}), and are based on the concept of metastable states and transfer operator (or transition matrix). Broadly speaking, $\Omega$ is decomposed into metastable sets, which are sets that are long-lived and in which the system gets trapped. Then a transition matrix $P_{ij}(\tau)$ is defined using the probability to reach a metastable set $j$ if one starts a trajectory of length $\tau$ (the lag-time) in set $i$. The analysis of the eigenvalues lead to the concept of Perron cluster. From the eigenvectors and eigenvalues, one can derive the rate and other kinetic information.

Although derived apparently independently and at a later date, some groups started exploring how one could model molecular systems using Markov state models, a well-known theory but which has been only (relatively) recently applied to modeling bio-molecular systems. See \cite{Singhal:2004is,Swope:2004gx,Swope:2004da,chodera2007automatic,chodera2006long,Noe2007}. See related work by~\cite{Shalloway:1996en}. Many of the theory and results for Markov state models can be found in the literature on conformation dynamics. In some sense, Markov state models can be viewed as a practical implementation of conformation dynamics, that attacks the high-dimensionality of $\Omega$ by subdividing the space into a ``small'' number of cells, also called macro-states.

Specifically, the transition matrix $P_{ij}(\tau)$ giving the probability of reaching macro-state $j$ when starting from state $i$ after some lag time $\tau$ is used to compute the rate. Its first eigenvalue is equal to 1 and corresponds to the equilibrium distribution. The second eigenvalue is very close to 1 and can be used to estimate the relaxation rate, as well as the forward (reactant to product) and backward (product to reactant) rates. An important issue is the effect of the lag time $\tau$. At short times, non-Markovian effects, or memory, are present, that is the Markov state model is not accurate and the estimated rate suffers from a systematic bias.

\medskip

{\bf Reactive trajectory sampling.} The last class of methods groups two separate approaches that in some sense combine ideas from transition path sampling and a subdivision of space similar to Markov state models. One such method, called weighted ensemble Brownian dynamics, originates in \cite{Huber:1996dn}. Although this paper is similar in spirit to transition interface sampling or milestoning, it can be easily extended to a general partitioning of space, using for example Voronoi cells. This is an important extension since, as a result, the method remains efficient in cases where multiple pathways contribute to the rate or when the most important pathway is not known. This approach leads to a sampling of transition pathways between $A$ and $B$ and therefore does not rely on the Markovian assumption made in Markov state models. However like Markov state models, the efficiency of the sampling is improved by partitioning space into macro-states. A large number of walkers (simulations) are run in each macro-state. In order to maintain the population of walkers in each macro-state, a procedure was created to kill walkers in macro-states that are too crowded, and to split walkers when the number of walkers becomes too low. This method was recently revisited by \cite{Zhang:2010kf, Bhatt:2010df} who showed how the original approach could be extended.

The technique of non-equilibrium umbrella sampling of \cite{Warmflash:2007dz, Dickson:2009gt} is similar in spirit. It applies, like weighted ensemble Brownian dynamics, to non-equilibrium systems and systems with memory (the Markovian approximation is not required) and uses a partitioning of $\Omega$ into macro-states. In each macro-state, a large number of walkers are simulated. Each time a walker attempts to leave a macro-state, its position (and velocity if needed) is recorded. Then in order to restart (continue) the simulation, a random position is chosen from the set of walkers who attempted to enter this macro-state from other macro-states. 

\medskip

In the second half of the paper we will discuss in more details weighted ensemble Brownian dynamics, renamed Reactive Trajectory Sampling (RTS) to reflect its broader application. The original method, \cite{Huber:1996dn}, involves a procedure to split and kill walkers. We will revisit this method and propose an optimal procedure which leads to walkers with identical weights in each macro-state, a strategy which minimizes the statistical errors. We will discuss how the choice of macro-states affects the statistical errors and what the optimal choice is. We will present a new ``multi-colored'' algorithm that allows computing multiple rates (eigenvalues of the transfer operator from conformation dynamics) and accelerate convergence. This approach shares some similarities with the technique of core sets, \cite{2011JChPh.134t4105S}.

We will propose a novel error analysis of Markov state models, by considering the sensitivity of the eigenvalues to perturbations in the transition matrix. This will lead to estimates of the systematic errors (non-Markovity) and statistical errors, and their dependence on the lag-time $\tau$ (length of trajectories used to calculate the transition matrix). The choice of macro-states influences the decay of the non-Markovian effects; we will discuss what the optimal choice is. This optimal choice is, as can be expected, difficult to realize in practice, but this provides nonetheless a guideline to improve and refine an initial guess. We will make an argument showing that in the general case, statistical errors increase with $\tau$, showing that an optimal tradeoff must be found between memory effects (small $\tau$) and statistical errors (large $\tau$). The reader is referred to \cite{sarich2010approximation,Prinz:2011id,Hinrichs:2007hj,Singhal:2005dk,2009PhRvE..80b1106M} for papers that discuss the numerical errors in Markov state models.

Some numerical results on simple benchmark problems in 1D and 2D are given at the end to illustrate the numerical properties of Markov state models and RTS.

\medskip

Since there is a large number of methods to choose from, with different strengths and weaknesses, we attempt to summarize their main features for the purpose of comparing these methods together. We considered three axes in our evaluation: generality, computational cost, and parallel scalability. Each characteristic is ranked low, medium, or high. Generality relates to the number of approximations or assumptions that are required by the method to be accurate. For example, the method of reactive flux focusses its analysis on the transition region, typically a saddle point. When the assumptions are satisfied the calculation may be very accurate. However, in some instances the prediction may not be satisfactory and need to be improved for example using a method based on transition path sampling. The computational cost should also be taken as a general guideline since it will vary tremendously depending on the system. However some overall conclusions can be made regarding computational cost. Typically the situation is that generality is traded for computational cost. Finally since all large scale calculations require a parallel computer (multicore processors, graphics processing units, and parallel clusters), we also rank methods according to the amount of parallelism they offer, although in all cases, the amount of parallelism is very large and scalability is typically not an issue. Codes that run on slow networks (grid computing, cloud computing) will be more sensitive to these issues. Finally, the information in the list below should be taken merely as a guideline since all conclusions are in general strongly system and implementation dependent.

\bigskip

\noindent {\bf Reactive flux} \\
\underline{Generality:} low. These approaches are typically the ones that require the greatest amount of knowledge about the system and relatively strict assumptions, in particular regarding the energy landscape near the transition region. \\
\underline{Computational cost:} low. A free energy calculation is required. Although this can be difficult in some instances, this type of calculation is typically easier than with the other methods. \\
\underline{Parallel scalability:} medium. Many methods are available to calculate free energy and are quite scalable. In most cases, it may become difficult to increase the number of processors compared to the other techniques described here.

\medskip

\noindent {\bf Transition path sampling} \\
\underline{Generality:} high. These approaches are among the most ``direct'' and require little or no assumption. \\
\underline{Computational cost:} high. As a result of being very general, they lead to extensive sampling and typically long simulation times. \\
\underline{Parallel scalability:} high. In most cases, one can sample paths independently.

\medskip

\noindent {\bf Transition interface sampling, forward flux sampling} \\
\underline{Generality:} high. These approaches also require little or not assumption. Forward flux sampling in addition only requires the ability to run forward simulations making them applicable in situations where transition interface sampling fails. \\
\underline{Computational cost:} high. A large number of pathways must be sampled. Forward flux sampling may in some cases converge a little slower, in particular when the sampling of the initial interfaces is poor or turns out to be insufficient for later interfaces. These methods typically perform better than transition path sampling. \\
\underline{Parallel scalability:} medium. Processing the interfaces is sequential, making the method less parallel. 

\medskip

\noindent {\bf Milestoning} \\
\underline{Generality:} medium. Strong assumptions are required regarding loss of memory between milestones. As discussed in the main text, the original method has been extended and made more general. \\
\underline{Computational cost:} medium. Only short pathways between adjacent milestones are required making this approach less expensive than the two previous methods. \\
\underline{Parallel scalability:} high. The milestones can be processed independently.

\medskip

\noindent {\bf Markov state models} \\
\underline{Generality:} medium. The Markovian assumption must apply, which depends on the choice of Markov (or macro) states and the lag time between observations. \\
\underline{Computational cost:} medium. The sampling is mostly local within each macro state, leading to efficient sampling. \\
\underline{Parallel scalability:} high. Macro states are processed independently and require running a large number of short trajectories (in some implementations).

\medskip

\noindent {\bf Weighted ensemble Brownian dynamics and non-equilibrium umbrella sampling.} \\
\underline{Generality:} high. The Markovian assumption is not required. The level of generality is similar to transition path sampling and related methods. The convergence of weighted ensemble Brownian dynamics is relatively easy to monitor. \\
\underline{Computational cost:} high. The fact that the Markovian assumption does not apply typically leads to higher computational cost compared to Markov state models. \\
\underline{Parallel scalability:} high. It is similar to Markov state models, although some communication is required to update the weight of walkers. This involves a small amount of communication, but one that occurs at regular intervals.

\bigskip

Weighted ensemble Brownian dynamics is in many respects similar to transition path sampling and differs ``mostly'' in the technique used to generate paths joining $A$ and $B$ and enhance the sampling. In this paper, the method of weighted ensemble Brownian dynamics will be referred to as reactive trajectory sampling to indicate that more general formulations have been created since the original paper of~\cite{Huber:1996dn}.

Many theoretical results for these methods have been proved in the context of Langevin dynamics or Brownian dynamics. Results are often derived in the context of Brownian dynamics (over-damped dynamics). Extensions to Langevin are in most cases possible, although the proofs become more technical. In some cases, a method may depend only on some stochastic process with minimal assumptions. The theoretical derivations at the end of this paper are done in the context of Brownian dynamics but results can be extended to Langevin dynamics. Extensions to Newtonian dynamics (deterministic) are much more difficult and in most cases these extensions do not exist yet.

\bigskip

\topcaption{Notations used in this paper. Not all the notations used in this paper can be found here. Notations that are local to a page or paragraph have been omitted.
\label{tab1} }
\begin{supertabular}{ll}
TIS & transition interface sampling \\
FFS & forward flux sampling \\
MSM & Markov state models \\
WEB & weighted ensemble Brownian dynamics \\
RTS & reactive trajectory sampling \\
MFEP & minimum free energy pathway \\	
$\Omega$ & conformational space \\
$x_k$ & coordinates of atom $k$ \\
$m_k$ & mass of atom $k$ \\
$U(x)$ & potential energy \\
$t$ & time \\
$T$ & temperature \\
$\beta$ & $\beta = (kT)^{-1}$ \\
$\langle \; \rangle$ & statistical average \\
$\rho(x)$ & $\rho(x) = e^{-\beta U(x)}/Z$ \\
$Z$ & partition function, $Z = \int e^{-\beta U(x)} \; dx$ \\
$A$ & subset of $\Omega$; reactant states \\	
$B$ & subset of $\Omega$; product states \\
$k_{AB}$ & rate from $A$ to $B$ \\
$\tau_{AB}$ & mean passage time from $A$ to $B$ \\
$k_\text{TST}$ & rate as predicted from transition state theory \\
$k_\text{Kramers}$ & rate as predicted from Kramers' theory \\
$\dot{f}(t)$ & time derivative of $f(t)$ \\
$\tau$ & Lag-time in Markov state models \\
$P_{ij}(\tau)$ & probability to be in state $j$ when starting from $i$ after time $\tau$ \\
$\xi$ & reaction coordinate or order parameter that monotonically increase from $A$ to $B$ \\
$(\xi_1,\ldots,\xi_p)$ & set of generalized coordinates \\
$A(\xi)$ & free energy \\
$D$ & diffusion tensor \\
$\chi_A$ & characteristic function of some set, $A$ \\
$C(t)$ & conditional probability to find the system in $B$ at $t$ provided it was in $A$ at time 0 \\
$x(\mathcal{T})$ & a discrete trajectory in $\Omega$ \\
$\mathcal{P}(x(\mathcal{T}))$ & probability density function of trajectories \\
$\mathcal{P}_{AB}[x(\mathcal{T})]$ & probability density function for the transition path ensemble \\
$S_i$ & hypersurface in $\Omega$; used in transition interface sampling, forward flux sampling, etc. \\
$V_i$ & cell or macro-state in Markov state models and related methods \\
& In many cases $S_i$ is the set of points such that $\xi(x) = \xi_i$ \\
$\Phi_{A,1}$ & flux from region $A$ to $S_1$ \\
$P_A(\xi_{i+1} | \xi_i)$ & for trajectories coming from $A$, probability to reach $S_{i+1}$ starting from $S_i$ \\
$\pi(x)$ & committor function \\
$\mu_i$ & eigenvalue of the transition matrix \\
$\lambda_i$ & $\lambda_i = -(\ln \mu_i)/\tau$; they are often an approximation of the eigenvalues \\
& of the Fokker-Planck equation \\
$\mathbf{P}(\;)$ & used to denote the probability of some event happening \\
$X_n$, $X_t$ & homogeneous Markov process indexed by $n$ or $t$ \\
$N_{ij}$ & number of observed crossings from macro-state $i$ to $j$ \\
$T_i$ & length of simulation in macro-state $i$ \\
$W_i$ & statistical weight of macro-state $i$ \\
$w_i$ & statistical weight of walkers in weighted ensemble Brownian dynamics \\ & or reactive trajectory sampling \\
$\rho(x,t|x_0,0)$ & probability to be at $x$ at time $t$ if the system was at $x_0$ at time 0 \\
$\rho_k$ & eigenfunction of the forward Fokker-Planck equation \\
$\psi_k$ & eigenfunction of the backward Fokker-Planck equation \\
$\lambda_k$ & corresponding eigenvalues \\
\end{supertabular}

\vspace*{2em}

\section{Transition path sampling}

\subsection{Reactive flux and transition state theory}

We start the discussion with the method of reactive flux which is a long standing approach to computing reaction rates. The idea goes back to \cite{Marcelin:1915uv, Eyring:1931wh, Eyring:1935dt, Horiuti:1938ti, Wigner:1938uw}, who developed the initial theory of chemical reaction kinetics. We outline the main ingredients in this type of approach. We assume that region $A$ is a subset of the conformational space of the molecular system and that it represents in the system in its reactant state. Similarly $B$ denotes the region defining the product states. Analytical approximation for the rate can be obtained if one assumes that a coordinate $\xi$ can be defined which describes the reaction, a reaction coordinate. It is assumed that when $\xi=0$ the system is in $A$ and when $\xi=1$ the system is in $B$. The value $\xi=\xi^*$ corresponds to the transition region or barrier between $A$ and $B$. 

We define the characteristic function $\chi_A$ (resp.\ $B$) which is 1 in the set $A$ and 0 outside. Then using these functions, we can express the conditional probability to find the system in state $B$ at time $t$ provided it was in $A$ at time 0:
\begin{equation}
	C(t) = \frac{\langle \chi_A[\xi(0)]\chi_B(\xi(t) \rangle}{\langle \chi_A \rangle} 
	\label{eq:C}
\end{equation}
Brackets $\langle \, \rangle$ are used to denote a statistical average. Regions $A$ and $B$ are separated by a transition region and the rate is determined by the rate at which this transition or barrier is crossed. At the molecular scale, there is some correlation time $\tau_\text{mol}$ associated with this crossing. That is for times larger than $\tau_\text{mol}$, the system has forgotten how it went from $A$ to $B$. Then for times $t$ between $\tau_\text{mol}$ and the reaction time $\tau_\text{rxn}$, $\tau_\text{mol} < t \ll \tau_\text{rxn}$, the time derivative of $C(t)$, called the reactive flux, reaches a plateau (\cite{Chandler:1978hh}), and
\begin{equation}
	\dot{C}(t) \approx k_{AB}
\end{equation}
The symbol \; $\dot{\,}$ \; denotes a time derivative.

Using transition state theory (TST), under the assumption that the recrossing of the barrier between $A$ and $B$ can be neglected, one can derive an expression for $k_{AB}$ using Eq.~\eqref{eq:C} (\cite{Chandler:1978hh,Chandler:1987tp}):
\begin{equation}
	k_\text{TST} = \frac{1}{2} \langle |\dot{\xi}| \rangle_{\xi=\xi^*}
	\frac{e^{-\beta A(\xi^*)}}{\int_{-\infty}^{\xi^*} e^{-\beta A(\xi)} d\xi}
\end{equation}
where $A(\xi)$ is the free energy, and $\langle \rangle_{\xi=\xi^*}$ denotes an ensemble average with $\xi$ constrained at $\xi^*$. This approach has some drawbacks. It always overestimates the rate. It requires a good reaction coordinate and a precise determination of the free energy maximum to locate the barrier. Nevertheless the method is computationally efficient and involves only quantities that can be computed with relatively low computational cost. Among the many methods to calculate the free energy (in this context the potential of mean force), see for example \cite{Lelievre:2010td, Anonymous:2007wt, 2001JChPh.115.9169D, 2004JChPh.121.2904H, 2004JChPh.120.3563R, 2007JCoPh.222..624L, 2008JChPh.128n4120D}.

Related approaches include Kramers' rate theory (\cite{Gardiner:1997tb,Hanggi:1990en}), which was developed in the context of Langevin equations and overdamped dynamics. There are many connections between transition state theory and Kramers' theory, \cite{Hanggi:1990en}. In particular Kramers' rate can be related to the ``simple'' TST rate through:
\begin{equation}
	k_\text{Kramers} = \frac{\lambda_+}{\omega_\text{bar}} \; k_\text{TST}
\end{equation}
In this expression the potential at the transition point is assumed to be locally quadratic with stiffness $\omega_\text{bar}^2 = -m^{-1} U''(x_\text{bar})$ ($m$ is the mass of the particle in a 1D model), and $\lambda_+$ is a function of the friction in the Langevin model and $\omega_\text{bar}$. It can be shown that $k_\text{Kramers}$ is equal to the multidimensional TST rate for a heat bath describing strict Ohmic friction, \cite{Hanggi:1990en}, pp. 268 \& 272. As the friction in the Langevin model goes to zero $\lambda_+ \to \omega_\text{bar}$ and $k_\text{Kramers} \to k_\text{TST}$. Moreover we always have $k_\text{Kramers} < k_\text{TST}$. The rate $k_\text{Kramers}$ is itself an upper bound on the true rate given by
\begin{equation}
	k(t) = \frac{\langle \dot{\xi}(0) \theta(\xi(t)-\xi^*)\rangle_{\xi(0)=\xi^*}}{\langle \theta(\xi^*-\xi(0))\rangle} 
\end{equation}
where $\theta$ is the Heaviside function.

This basic approach using TST has been improved in many ways including the use of harmonic approximations to model the minimum energy basins and transition region; see \cite{Dellago2009, Gardiner:1997tb, Hanggi:1990en}. In variational TST, one attempts to improve the predicted rate by finding a dividing surface between $A$ and $B$ that minimizes the rate, see for example \cite{Truhlar:1984vh, Tucker:1995us}. The overestimation of the rate by TST is a result of neglecting the re-crossing of the dividing surface. Several authors have proposed corrections to the basic TST approach to account for these effects. See for example \cite{Bennett:1977ui, Chandler:1978hh}.

\subsection{Transition path sampling}

The method of reactive flux is attractive as its computational cost is often tractable. The primary calculation is obtaining the free energy profile along the reaction coordinate $\xi$. The accuracy is very dependent on the choice of coordinate. Specifically, the transmission coefficient, which measures the amount of re-crossing taking place, is dependent on the definition of $\xi$ and a low transmission coefficient leads to inaccuracies or inefficiencies. Transition path sampling were in part proposed to alleviate the need to define this coordinate as they rely primarily on sampling trajectories going from $A$ to $B$ with no knowledge of $\xi$ required. However we will see that later versions, again, are dependent on some knowledge of $\xi$.

Transition path sampling (\cite{Bolhuis:2002ew,1998JChPh.108.1964D,Dellago:2002uf}) is a Monte-Carlo method that allows sampling the ensemble of trajectories. For each discrete trajectory
\begin{equation}
	x(\mathcal{T}) = \{ x_0, x_{\Delta t}, \ldots, x_{\mathcal{T}} \}
\end{equation}
one can define a probability $\mathcal{P}(x(\mathcal{T}))$ to observe such a trajectory. Then the transition path ensemble defines a probability density in the space of trajectories that is non-zero only for trajectories that connect $A$ and $B$, and therefore its probability density function is defined as:
\begin{equation}
	\mathcal{P}_{AB}[x(\mathcal{T})] = \chi_A(x_0) \chi_B(x_{\mathcal{T}}) \mathcal{P}(x(\mathcal{T})) / Z_{AB}({\mathcal{T}})
\end{equation}
where $Z_{AB}$ is the appropriate normalization factor.

\begin{figure}[htbp]
\centering
\begin{tikzpicture}
\draw (-3,0) ellipse (2cm and 1cm) node {$A$};
\draw (3,0) ellipse (2cm and 1cm) node {$B$}; 
\draw [thick,red,->] (-2.5, 0.1) .. controls +(10:1) and +(160:1) .. (0,1)
	.. controls +(-20:0.5) and +(190:1) .. (1,0.5)
	.. controls +(10:1) and +(110:1) .. (2,-0.4);
\draw [thick,blue,->] (-2.5, -0.3) .. controls +(-20:1) and +(200:1) .. (-.5,-0.2)
	.. controls +(20:0.5) and +(210:0.5) .. (1,0.5)
	.. controls +(30:1) and +(120:1) .. (3.5,0.2);	
\filldraw (1,0.5) circle (0.08) node [text width = 1.8in,anchor=south,yshift=0.2in]{Point at which perturbation in momentum is applied};	
\end{tikzpicture}
\caption{Shooting method of transition path sampling. The red curve is the old trajectory. The blue curve is obtained by perturbing the momentum at some point and regenerating a trajectory by using a forward and backward time integrator. \label{fig:shoot}}
\end{figure}
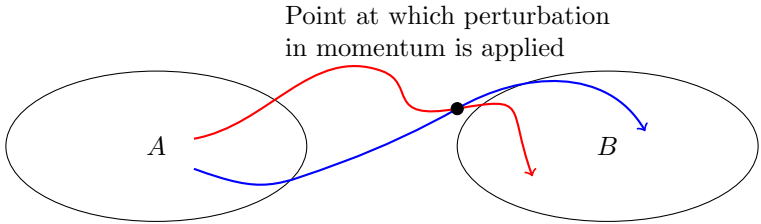

Trajectories in this ensemble can be generated using a Monte-Carlo procedure. A popular method is the so-called shooting method, in which a point is selected along the trajectory, the momentum is modified and a new trajectory is generated using a forward and backward time integration. If a stochastic dynamics is used then no perturbation is necessary since the random number generator will lead to a different trajectory. Then the appropriate acceptance probability is used (Metropolis-Hastings algorithm) to accept or reject this new trajectory. See Fig.\ref{fig:shoot}. Trajectories that do not start in $A$ and end up in $B$ are always rejected. This procedure can be improved by adding extra moves like shifting moves or path reversal moves, \cite{Dellago:2002uf}.

Many thermodynamic and kinetic properties can be determined from transition path sampling. As before the rate is related to $C(t)$ through $k_{AB}(t) = \dot{C}(t)$ and (\cite{Dellago2009, Dellago:2002uf})
\begin{equation}
	C(t) = \frac{\int \mathcal{D}x(t) \mathcal{P}(x(t)) \chi_A(x_0) \chi_B(x_t)}{\int \mathcal{D}x(t) \mathcal{P}(x(t)) \chi_A(x_0)} 
\end{equation}
where $\int \mathcal{D}x(t)$ is an integration over all possible paths (precisely, all possible points along the approximation at discrete time steps of the continuous trajectory).

There is a numerical difficulty in computing
\begin{equation}
	\int \mathcal{D}x(t) \, \mathcal{P}(x(t)) \, \chi_A(x_0) \, \chi_B(x_t)
\end{equation}
since for most trajectories of length $t$, $\chi_A(x_0)=0$ or $\chi_B(x_t)=0$. A common procedure to make this tractable is to use umbrella sampling and an order parameter $\xi$ which approximates the reaction coordinate. Denote $P_A(\xi',t)$ the probability that a trajectory started in $A$ is such that $\xi=\xi'$ at time $t$, then (\cite{Dellago2009}):
\begin{equation}
	C(t) = \int_{\xi_\text{min}^B}^{\xi_\text{max}^B} P_A(\xi,t) \, d\xi
\end{equation}
where $\xi_\text{min}^B \le \xi \le \xi_\text{max}^B$ defines region $B$.
The quantity $P_A(\xi,t)$ will be difficult to obtain for values of $\xi$ that are close to 1/2 (barrier) or beyond. The method of umbrella sampling (\cite{Torrie:1977wa}) can then be conveniently used. For a window $W_i$ in the interval $[0,1]$ define:
\begin{equation}                             
	P_A^{W_i}(\xi,t) \define
	\frac{\displaystyle \int \mathcal{D}x(t) \, \mathcal{P}(x(t)) \, \chi_A(x_0) \, \chi_{W_i}(x_t) \, \delta[\xi-\xi(x_t)]}
	{\displaystyle \int \mathcal{D}x(t) \, \mathcal{P}(x(t)) \, \chi_A(x_0) \, \chi_{W_i}(x_t)}
\end{equation}
Then each $P_A^{W_i}(\xi,t)$ can be efficiently computed for small enough windows, while $P_A(\xi,t)$ over the entire $[0,1]$ interval is obtained by patching together the different $P_A^{W_i}(\xi,t)$ and matching the curves to produce a single profile.

\medskip

{\bf Transition path theory.} Transition path theory (TPT) is a theoretical framework to study transition path ensembles. It considers so-called reactive trajectories, which are the trajectories sampled by transition-path sampling. The theory is derived in the context of Langevin and over-damped dynamics. TPT provides definitions for the probability density of reactive trajectories. Most results involve the forward $q_+$ or backward $q_-$ committor functions. The forward committor function for example is defined as the probability that, starting from some point outside of $A \cup B$, the system first reaches $B$ before $A$, \cite{E:tb}. As an example, the probability density of reactive trajectories is defined as:
\begin{equation*}
	\rho_\text{reactive}(x) = \rho(x) \, q_+(x) \, q_-(x)
\end{equation*}
where $\rho(x)$ is the equilibrium distribution for the process. More importantly expressions are provided for the probability current of reactive trajectories. This allows identifying important transition pathways and transition tubes. The example of the maze is helpful. The path in the maze with the largest reactive flux corresponds to the shortest path from $A$ to $B$. Actual transition pathways take many detours and visit many dead ends before back-tracking and ending up in $B$. Analyzing the reactive flux allows bypassing these detours and take a straight route to $B$. In practice this allows identifying important reaction mechanisms joining $A$ and $B$. Expressions are also given for the reaction rate, basically by integrating the probability flux over a subdividing surface.

\cite{Metzner2009} presents an application to discrete Markov processes. In this case the discrete probability current for reactive trajectories is given by:
\begin{equation*}
	f_{ij}^{A \to B} =
	\begin{cases}
		\rho_i \; q^-_i \; l_{ij} \; q^+_j, & \text{if $i \neq j$,} \\
		0, & \text{otherwise.}
	\end{cases}
\end{equation*}
where $l_{ij}$ is the infinitesimal generator (rate matrix) of a continuous-time Markov chain.

See~\cite{E:2006hg} for an earlier work. A series of illustrative examples are proposed in~\cite{Metzner:2006du}.

\subsection{Transition interface sampling}

One issue in the previous approach is the fact that the trajectories need to have a fixed length $t$. See~\cite{Dellago:2002uf} (Section IV D ``A Convenient Factorization'') for an algorithm that relaxes this requirement. This length must be chosen carefully. A short time will lead to inaccuracies while a long time leads to a larger computational cost to update the trajectories. The method of transition interface sampling (TIS) introduces a novel concept, somewhat related to the previous algorithm with umbrella sampling, in which $n+1$ multi-dimensional surface in the space outside of $A$ and $B$ are defined. These surface are such that they can be used roughly to measure the progress of the reaction. For example, for some order parameter $\xi$ and choosing a value $\xi_i$, we can define a surface $S_i$ by the equation $\xi(x) = \xi_i$. We assume that $\xi_0 < \xi_1 < \ldots < \xi_n$; $S_0$ is the boundary of $A$ and $S_n$ is the boundary of $B$. Several methods (transition interface sampling and forward flux sampling) start from a reinterpretation of the rate as:
\begin{equation}
	k_{AB} = \Phi_{A,1} \; P_A(\xi_n | \xi_1)
\end{equation}
The first term, $\Phi_{A,1}$, is the flux of trajectories going out of $A$ and crossing the first surface $S_1$. This quantity can be obtained by running a long trajectory (ignoring parts of the trajectory that last visited $B$) and counting the number of times the surface $S_1$ is crossed, per unit time. Only positive crossing are counted (that is moving away from $A$ towards $B$) and re-crossings are ignored until $A$ is entered again (that is crossings for which the previous crossing was $S_1$, not $A$, are ignored). The second term, $P_A(\xi_n | \xi_1)$ is the probability to reach $S_n$, assuming that the trajectory has crossed $S_1$ and that it does not cross $A$ before $S_n$. See \cite{vanErp:2003ik}.

The probability $P_A(\xi_n | \xi_1)$ is difficult to calculate since it is very small. This can be remedied using the intermediate surfaces $S_i$, $2 \le i \le n-1$. By construction, the function $\xi$ being continuous and since $\xi_0 < \xi_1 < \ldots < \xi_n$, it is not possible to last come from $A$ and cross $\xi_i$ without having crossed first $\xi_{i-1}$. Using this result, it is possible to show that (\cite{vanErp:2003ik}):
\begin{equation}
	P_A(\xi_n | \xi_1) = \prod_{i=1}^{n-1} P_A(\xi_{i+1} | \xi_i)
\end{equation}
where $P_A(\xi_{i+1} | \xi_i)$ is the probability to cross $S_{i+1}$ before $A$ assuming the trajectory had previously crossed $A$ then $S_i$. Although $P_A(\xi_n | \xi_1)$ can be very small, the quantities $P_A(\xi_{i+1} | \xi_i)$ are much larger and can be reliably obtained through direct sampling. The procedure to calculate $P_A(\xi_{i+1} | \xi_i)$ is similar to the procedure for TPS. It is illustrated in Fig.~\ref{fig:tis}.

\subsection{Forward flux sampling}

The forward flux sampling method (FFS) was conceived by \cite{Allen:2005dn, Allen:2006ch, Allen:2006cp}. An earlier paper by \cite{Harvey:1993de} bears some conceptual similarities. In TIS, one generates new trajectories by perturbing a point and then integrating forward and backward. In the forward flux sampling method, only forward integration is used. This can be essential in cases where backward integration is not possible. For example in non-equilibrium systems, the lack of detailed balance and absence of time-reversal symmetry means that TIS (or milestoning which will be described in the next section) is not applicable. FFS is one of the few methods applicable to such systems.

FFS uses the same basic framework as TIS expressing the rate as
\begin{equation}
	k_{AB} = \Phi_{A,1} \; \prod_{i=1}^{n-1} P_A(\xi_{i+1} | \xi_i)
\end{equation}
FFS start by calculating $P_A(\xi_2 | \xi_1)$ using trajectories initiated from $S_1$. Some of these trajectories may fail to reach $S_2$ (i.e., reach $A$ before $S_2$) while others will reach $S_2$ successfully (before reaching $A$). For those that reach $S_2$, the first hitting point (first point where the trajectory crosses $S_2$) is saved. Those points are then used to calculate the next conditional probability, $P_A(\xi_3 | \xi_2)$. FFS therefore only requires forward integration of trajectories. The starting points on $S_i$ are in some sense fixed and are produced solely as a result of the sampling during the calculation of $P_A(\xi_i | \xi_{i-1})$. Typically the computation for $P_A(\xi_i | \xi_{i-1})$ is continued until satisfactory accuracy is achieved and enough points have been generated on $S_i$. We note that the forward flux approach not only yields the rate constant, but also the complete transition pathways, which can be reconstructed by ``gluing'' the successful trajectory segments together. The method is depicted in Fig.~\ref{fig:ffs}.

A limitation of this method is that the accuracy at later interfaces depends on the sampling at earlier interfaces. For example if the first interface is relatively poorly sampled the error will propagate throughout the next interfaces. In addition it is very possible that initial trajectories with low probabilities, in the end, make large contribution to the flux. Such a situation would lead to a large standard deviation and statistical errors. 

Example applications are discussed for example in \cite{Valeriani:2007hv,Borrero:2007eq,Allen:2005dn}. \cite{Allen:2006cp} describe FFS along with two other methods, the branched growth method and the Rosenbluth method. \cite{Allen:2006ch} proposes an analysis of the efficiency of these methods.

\subsection{Milestoning}

Similar to the previous method, the milestoning technique of \cite{Faradjian:2004cr,West:2007dg} is based on a set of separating hyper-surfaces that are used to measure the progress of the system from $A$ to $B$. The advantage of this method is that trajectories need only to be run from one interface $S_i$ (or milestone) to the next $S_{i+1}$ or previous interface $S_{i-1}$. In contrast, the previous methods required running trajectories until they reach the next milestone or set $A$. Milestoning is inherently more parallel since each interface $S_i$ can be processed in parallel with the others, whereas the other methods require a sequential processing of the interfaces. A drawback is that the interfaces need to be sufficiently separated so that we can assume that the system loses memory in the time it takes to reach the next interface.

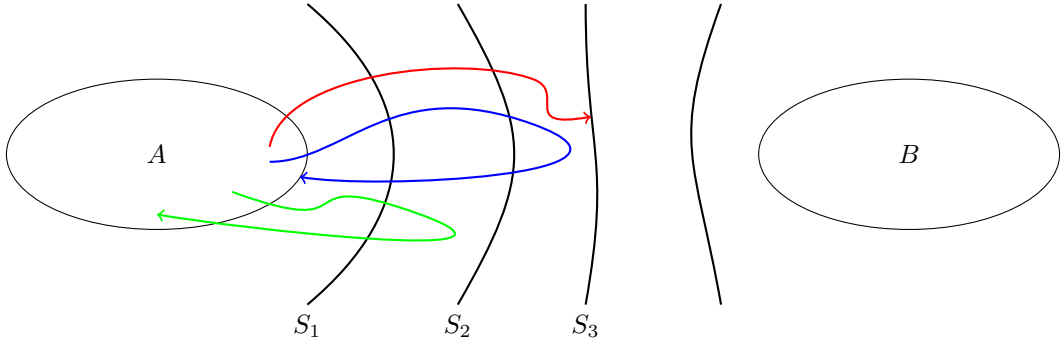
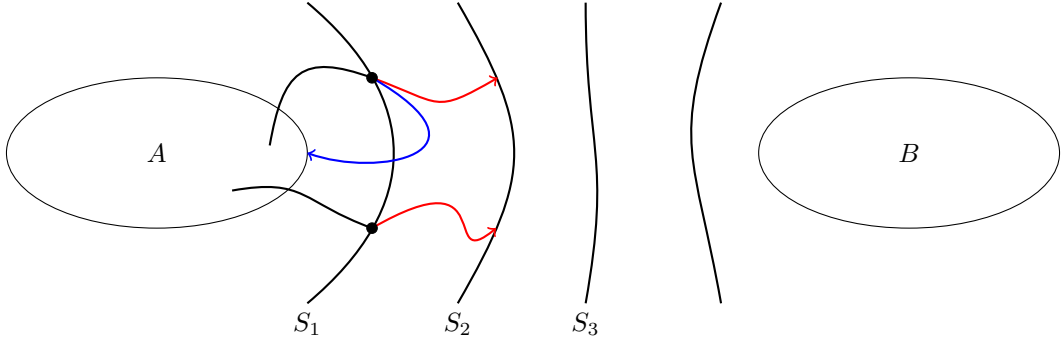
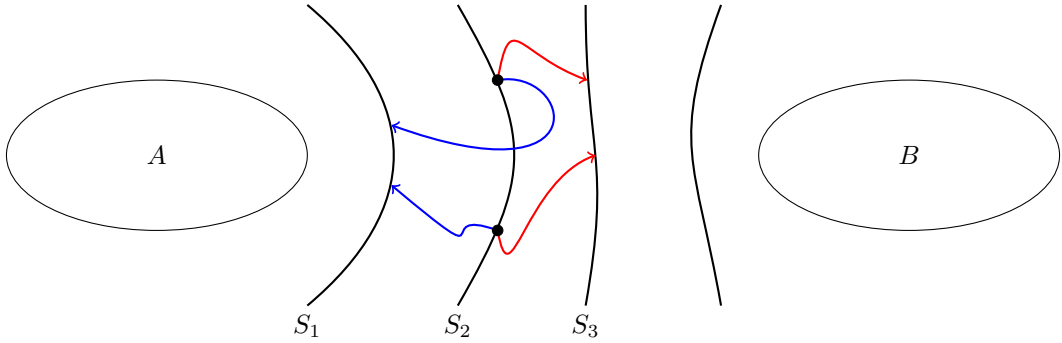
\begin{figure}
\centering

\subfloat[Transition interface sampling algorithm. Calculation of $P_A(\xi_3|\xi_2)$. The sampling is similar to TPS. The main difference is that the only paths that are considered are those that last come from $A$, cross $S_2$ and then cross $A$ or $S_3$. A new path is generated by perturbing a given path and then integrating forward or backward. The integration can be stopped as soon as $A$ or $S_3$ are reached. The new path is then accepted only if there is a segment crossing $S_2$ that last come from $A$. Green trajectory: it has not reached the surface $S_2$ and therefore is discarded. Red trajectory: it has reached $S_3$ before $A$ and therefore counts as 1. Blue trajectory: it has reached $A$ before $S_3$ and therefore counts as 0. \label{fig:tis}]{
\begin{tikzpicture}
\draw (-5,0) ellipse (2cm and 1cm) node {$A$};
\draw (5,0) ellipse (2cm and 1cm) node {$B$};
\draw [thick] (-3,-2) node [anchor=north] {$S_1$} .. controls +(40:2) and +(-40:2) .. (-3,2);
\draw [thick] (-1,-2) node [anchor=north] {$S_2$} .. controls +(60:2) and +(-60:2) .. (-1,2);
\draw [thick] (0.7,-2) node [anchor=north] {$S_3$} .. controls +(80:2) and +(-90:2) .. (0.7,2);
\draw [thick] (2.5,-2) .. controls +(100:2) and +(-110:2) .. (2.5,2);
\draw [thick,red,->] (-3.5, 0.1) .. controls +(80:1) and +(160:1) .. (0,1)
	.. controls +(-20:0.5) and +(190:1) .. (0.77,0.5);
\draw [thick,blue,->] (-3.5,-0.1) .. controls +(0:1) and +(160:2) .. (0,0.4)
	.. controls +(-20:2) and +(-10:1) .. (-3.1,-0.3);
\draw [thick,green,->] (-4,-0.5) .. controls +(-20:2) and +(160:2) .. (-1.5,-0.8)
	.. controls +(-20:2) and +(-10:1) .. (-5,-0.8);         
\end{tikzpicture}	 
}

\subfloat[Forward Flux Sampling. The black trajectories are used to generate starting points on $S_1$. The red trajectories are trajectories started from $S_1$ who made it to $S_2$, while the blue trajectories reached $A$ first. The red trajectories are used to generate starting points on $S_2$. \label{fig:ffs}]{	
\begin{tikzpicture}
\draw (-5,0) ellipse (2cm and 1cm) node {$A$};
\draw (5,0) ellipse (2cm and 1cm) node {$B$};
\draw [thick] (-3,-2) node [anchor=north] {$S_1$} .. controls +(40:2) and +(-40:2) .. (-3,2);
\draw [thick] (-1,-2) node [anchor=north] {$S_2$} .. controls +(60:2) and +(-60:2) .. (-1,2);
\draw [thick] (0.7,-2) node [anchor=north] {$S_3$} .. controls +(80:2) and +(-90:2) .. (0.7,2);
\draw [thick] (2.5,-2) .. controls +(100:2) and +(-110:2) .. (2.5,2);

\draw [thick] (-3.5, 0.1) .. controls +(80:1) and +(160:1) .. (-2.14,1);         
\draw [thick,red,->] (-2.14,1) .. controls +(-20:1) and +(210:1) .. (-0.47,1);         
\draw [thick,blue,->] (-2.14,1) .. controls +(-30:2) and +(-20:1) .. (-3,0);
\filldraw (-2.14,1) circle (0.07);

\draw [thick] (-4, -0.5) .. controls +(10:1) and +(160:1) .. (-2.14,-1);         
\draw [thick,red,->] (-2.14,-1) .. controls +(30:2) and +(220:1) .. (-0.49,-1);         
\filldraw (-2.14,-1) circle (0.07);
\end{tikzpicture}	 
}

\subfloat[Milestoning algorithm. Trajectories are initiated from $S_2$ using the equilibrium Boltzmann distributions. The exit times are recorded when the system reaches $S_3$ (red trajectories, function $K^+_2(t)$) or $S_1$ (blue trajectories, function $K^-_2(t)$). \label{fig:milestoning}]{
\begin{tikzpicture}
\draw (-5,0) ellipse (2cm and 1cm) node {$A$};
\draw (5,0) ellipse (2cm and 1cm) node {$B$};
\draw [thick] (-3,-2) node [anchor=north] {$S_1$} .. controls +(40:2) and +(-40:2) .. (-3,2);
\draw [thick] (-1,-2) node [anchor=north] {$S_2$} .. controls +(60:2) and +(-60:2) .. (-1,2);
\draw [thick] (0.7,-2) node [anchor=north] {$S_3$} .. controls +(80:2) and +(-90:2) .. (0.7,2);
\draw [thick] (2.5,-2) .. controls +(100:2) and +(-110:2) .. (2.5,2);

\draw [thick,->,red] (-0.47,1) .. controls +(80:1) and +(160:1) .. (0.72,1);
\draw [thick,->,blue] (-0.47,1) .. controls +(10:1) and +(-20:3) .. (-1.88,0.4);         
\filldraw (-0.47,1) circle (0.07);

\draw [thick,->,red] (-0.47,-1) .. controls +(-80:1) and +(200:1) .. (0.83,0);
\draw [thick,->,blue] (-0.47,-1) .. controls +(160:1) and +(-40:2) .. (-1.88,-0.4);        
\filldraw (-0.47,-1) circle (0.07);
\end{tikzpicture}	 
}
\caption{These schematic figures illustrate three different schemes: transition interface sampling, forward flux sampling, and milestoning.}
\end{figure}

The milestoning method proceeds by initiating trajectories at $S_i$, using the equilibrium distribution. Then it records the time it takes to reach $S_{i-1}$ or $S_{i+1}$. This leads to two time distribution functions $K_i^+(t)$ (to reach $S_{i+1}$) and $K_i^-(t)$ (to reach $S_{i-1}$). We define $K_i(t) = K_i^+(t) + K_i^-(t)$. See Fig.~\ref{fig:milestoning}. Then we can calculate two functions: $P_s(t)$, which is the probability of being at milestone $s$ (that is the last interface that was crossed was $S_s$), and $Q_s(t)$, the probability to transition to milestone $s$ at time $t$ (cross $S_s$ at $t$). Then these two functions satisfy the following set of integral equations:
\begin{align}
	P_s(t) & = \int_0^t \Big[
	1 - \int_0^{t-t'} K_s(\tau) \, d\tau
	\Big] Q_s(t') \; dt' \\
	Q_s(t) & = \eta_s \, \delta(t-0^+)
	+ \int_0^t \big[
	K^-_{s+1}(t-t') Q_{s+1}(t')
	+ K^+_{s-1}(t-t') Q_{s-1}(t')
	\big] \; dt'
\end{align}
In these equations $\eta_s$ is the initial milestone probability distribution. The first equation is saying that in order to be at $s$ the system must first transition to $s$ [$Q_s(t')$] and then stay at $s$ until time $t>t'$. The second equation says that to reach $s$ one must first get to $s-1$ and then transition from $s-1$ to $s$ (and similarly with $s+1$). The first term $\eta_s \, \delta(t-0^+)$ accounts for the initial distribution at $t=0$ ($\delta$ is the Dirac delta function). From these equations, the free energy and reaction rate can be obtained.

This approach is accurate once we assume that the system loses memory between milestones. With this, it becomes justified to independently generate initial conditions on each milestone following the Boltzmann distribution. These assumptions allow running independent calculations at each milestone and make it possible to run trajectories that stop as soon as the previous or next milestone is reached.

See a discussion of this approach in \cite{Elber:2005hm} as well as an application example to an allosteric transition with {\it deoxy Scapharca} hemoglobin in \cite{Elber:2007cc}.

\subsection{Milestoning using optimal milestones}

\label{milestoning}

The method of milestoning was recently revisited, see e.g., \cite{VandenEijnden:2008bn}. In this paper, it is shown that the assumption that the milestones need to be sufficiently far apart is not required provided that the interfaces $S_i$ are iso-surfaces of the committor function. The committor function, denoted $\pi(x)$, is the probability to reach $B$ before $A$ starting from $x$. In \cite{VandenEijnden:2009dw}, the milestoning method is extended to the case of a tessellation of the conformational space using Voronoi cells $V_i$. See Fig.~\ref{vi}. See a more complete discussion about Voronoi cells in Section~\ref{sec:msm}. In a more recent paper, \cite{Majek:2010kq} modify the original milestoning method using an approach that does not require a reaction coordinate (in that sense similar to \cite{VandenEijnden:2009dw}). However, instead of considering the committor function to guarantee the accuracy of the method, this approach focusses in guaranteeing a minimum separation between the milestones. This leads to greater accuracy since the assumption that memory is lost between crossing milestones is automatically satisfied.

We now review some of the mathematical underpinnings discussed in \cite{VandenEijnden:2008bn}. In particular we recall the main arguments to establish that the use of the committor function $\pi(x)$ to define the milestones leads to an exact rate prediction with milestoning, irrespective of whether the milestones are close or well separated. For this reason these milestones are called optimal milestones.

In \cite{VandenEijnden:2008bn}, one is concerned with computing the mean passage time $\tau_{AB}$ from $A$ to $B$. 

Consider for example the case of over-damped dynamics and the associated Fokker-Planck equation. The eigenvalues of the Fokker-Planck equations are denoted $\lambda_i$. We assume the $\lambda_i$ are ranked by magnitude so that $\lambda_1 = 0 < \lambda_2 < \lambda_3 < \ldots$. If there is a single eigenvalue $\lambda_2$ close to 0 ($\lambda_2 \ll \lambda_3$), then the mean passage time is related to $\lambda_2$ through $\tau_{AB} = \lambda_2^{-1} \rho(B)$ (where $\rho(B) = \int_B \rho(x)dx$).

Milestones, denoted by $S_i$, are defined as hypersurfaces in $\Omega$. The mean time to go from $A$ to $B$ can be obtained from the mean times to go from a milestone to another milestone. This has important consequences in terms of which assumptions need to be made to get an exact rate (in the absence of statistical errors). 

The key property is the following one. Let us assume we consider all the trajectories that go from milestone $S_i$ to $S_j$ ($j=i+1$ or $j=i-1$). Specifically we initialize trajectories on milestone $S_i$ with density:
\begin{equation}
	\frac{|\nabla \pi(x)| \; e^{-\beta U(x)}}{\int_{S_i} |\nabla \pi(x)| \; e^{-\beta U(x)} \; d\sigma_i(x)} 
	\label{fhpd}
\end{equation}     
where $d\sigma_i(x)$ is the surface element on milestone $S_i$. Then we can follow all the trajectories and record where they hit milestone $S_j$ (we discard trajectories that hit other milestones first). The density of points on $S_j$ is called the density of first hitting points.

In~\cite{VandenEijnden:2008bn}, it was shown that, if the milestones $\{S_k\}$ are iso-surfaces of $\pi$ (optimal milestones), then the density on $S_j$ is given by:
\begin{equation}
	\frac{|\nabla \pi(x)| \; e^{-\beta U(x)}}{\int_{S_j} |\nabla \pi(x)| \; e^{-\beta U(x)} \; d\sigma_j(x)}
\end{equation}
In particular this density is independent of the fact that the trajectories were started from $S_i$. This can be proved using the definition of the committor function and the forward Chapman-Kolmogorov equation.

The implication is that if we want to calculate the mean time to go from $i$ to $j$ we do not need to consider where the system is coming from. It is sufficient to initiate trajectories on $S_i$ with the density given above and calculate the average time required to reach $S_j$ (assuming this is the next milestone crossed).

We now discuss the calculation of the mean passage times and show that the property above is essential to derive expressions for the exact rate. Let us consider an absorbing boundary condition on the milestone, denoted $S_n$, that surrounds region $B$ (cemetery milestone~\cite{VandenEijnden:2008bn}). We assume that trajectories are initialized on $S_i$ with a probability density given by Eq.~\eqref{fhpd}. We define $T_i$ the mean length of a trajectory that start from $S_i$ and goes to $S_n$ (mean exit time) . We denote: $p_{ij}$ the probability that $S_j$ is crossed after $S_i$, and $\tau_i$ the mean time to hit any other milestone $S_j$, $j \neq i$. Using the result above regarding the first hitting point density, we have:
\begin{equation}
	T_i = \tau_i + \sum_{j \neq i, j \neq n} p_{ij} \, T_j
\end{equation}
This is true because the mean escape time for some $S_j$, $j \neq i$, is the same irrespective of whether the system comes from $S_i$ or not. \label{MSTproof} This is a weaker property than saying that the times to go from $S_i$ to $S_j$, $t_{ij}$, and $S_j$ to $S_k$, $t_{jk}$, are independent. One can construct examples where $t_{ij}$ and $t_{jk}$ are strongly correlated, while the density of first hitting points at a milestone is independent of the last milestone that was crossed.  

This equation can be written in matrix form as:
\begin{equation}
	(I - p^n) \; T^n = \tau^n
\end{equation}  
where $p_{ii}^n = 0$ and $p^n$ does not have the row or column
corresponding to $S_n$.

This equation requires computing $p_{ij}$ and $\tau_i$ for all milestones. This can be further simplified by observing that $p_{ij}/\tau_i$ has a simple interpretation. A formal proof can be given but here we simply outline the main points. We consider a very long trajectory. We will not discuss this further but limits must be taken as the trajectory length goes to infinity. Then:
\begin{gather}
	p_{ij} = \frac{\text{Number of times the system crossed $S_j$ after $S_i$}}
	{\text{Number of times the system crossed any milestone after last crossing $S_i$}} \\[1em]
	\tau_i = \frac{\text{Total length of time during which the last milestone crossed was $S_i$}}
	{\text{Number of times the system crossed any milestone after last crossing $S_i$}}
\end{gather}
So that:
\begin{equation}
	\frac{p_{ij}}{\tau_i} = \frac{\text{Number of times the system crossed $S_j$ after $S_i$}}{\text{Total length of time during which the last milestone crossed was $S_i$}}  
	\label{eq13}
\end{equation}
Assume that we now use a discrete integrator to integrate the dynamics (which may be Langevin or overdamped dynamics), with time step $\Delta t$. We define:
\begin{gather}
	P_{ij} = \mathbf{P}\text{(system crosses $S_j$ during the next step assuming it last crossed $S_i$)}
\\
P_{ii}=\mathbf{P}\text{(system crosses no
milestones during the next step assuming it last crossed
$S_{i}$)}
\end{gather}
Then from Eq.~\eqref{eq13}:
\begin{equation}
	P_{ij} = \Delta t \; \frac{p_{ij}}{\tau_i}
\end{equation}
if $\Delta t$ is small compared to the time required to go from a milestone to the next. This shows that:\footnote{The equation $\text{diag}(\tau)^{-1}(I-p) = \Delta t^{-1}(I-P)$ is also true along the diagonal since $p_{ii} = 0$ and $\Delta t/\tau_i = 1-P_{ii}$.}
\begin{equation}
	(P-I) \; \Big( \Delta t^{-1} T \Big) = -\one
	\label{mpt}
\end{equation}
The quantity $\Delta t^{-1} T_i$ is the mean number of steps required to go from $S_i$ to $S_n$.

The advantage of this equation is that it requires computing $P_{ij}$ only, which can be relatively easily computed. Assume that we have defined a partition of the conformational space $\Omega$ into cells $V_i$ such that the milestones $S_i$ form the boundary of these cells. See Fig.~\ref{vi}. Then the approach requires simply running independent simulations in all the cells $V_i$. For this, we need to use boundary conditions such that the system remains in the cell it started in, during the simulation. Let us assume that the trajectory hits a cell boundary with velocity $\dot{x}$. Then we know from the equilibrium probability density that there is another trajectory in the past or future, with a one-to-one mapping, which re-enters the cell through the same point. Its velocity can be chosen equal to $-\dot{x}$ (using the fact that the equilibrium probability density is even with respect to the momenta), or obtained using a hard wall reflection (now using the fact that the reflection conserves the equilibrium probability density) with
\begin{equation}
	\dot{x}_\text{re-entering particle}
	= \dot{x} - 2 \, (\dot{x} \cdot \boldsymbol{n}) \; \boldsymbol{n}
	\qquad
	\text{and $\boldsymbol{n}$ is the normal to the interface.}
	\label{eq:reflection}
\end{equation}
With this approach, one can generate a large number of samples in each cell, from which we can estimate $P_{ij}$:
\begin{equation}
	P_{ij} = \frac{\rho(V_a) N_{ij} / n_a}{\rho(V_a) N_i^a / n_a + \rho(V_b) N_i^b / n_b} 
	\label{eq14}
\end{equation}
where $N_{ij}$ is the number of times the system was found to cross $S_j$ after $S_i$, $a$ is the cell bordered by $S_i$ and $S_j$, while $b$ is the cell on the other side of $S_i$, $N_i^a$ (resp.\ $N_i^b$) is the number of steps for which the last milestone crossed was $S_i$ in cell $V_a$ (resp.\ $V_b$), and $n_a$, $n_b$ are the number of steps computed in each cell. This is basically a direct calculation of Eq.~\eqref{eq13}. Note that by construction, the density of first hitting points on the milestone is the exact one so that, up to statistical errors, Eq.~\eqref{eq14} is exact.

The advantage of this approach is its efficiency and the fact that the accuracy is more or less independent of the energy barrier between $A$ and $B$. It does not have any systematic error unlike the previous approaches. The main drawback is the requirement that the optimal milestones are iso-surfaces of the committor function, which again is difficult to realize in practice. 

This requirement was relaxed in \cite{Majek:2010kq}.

\section{Conformation dynamics and Markov state models}

\subsection{Conformation dynamics}

All the methods discussed above attempt in a sense to do a direct calculation of the rate, either by computing a reactive flux at the transition barrier or by sampling reactive trajectories going from $A$ to $B$. We now discuss another class of methods that also attempt to calculate the free energy and the rate but, indirectly, by calculating the rate of transition between metastable basins. If one derives a statistical model of the system in terms of hops or transition between states, then an eigenvalue analysis can be used to calculate reaction rates, metastable states, and extract many relevant kinetic and thermodynamic information.

The idea goes back to Deuflhard and Sch\"utte, who realized that computing time averages of physical observables or relaxation times of conformational changes (using molecular dynamics for example) was largely determined by the existence and properties of invariant sets, called metastable sets. These are by definition sets (subsets of $\Omega$, the conformational space of the molecular system) such that the system stays trapped in these sets for extended periods of times and with very rare transitions between sets. This has led to the conformation dynamics approach which aims at identifying these sets, and computing the transition rates between these sets. The first paper goes back to \cite{Deuflhard1998}, although the term ``essential dynamics'' can be found in~\cite{Amadei1993}; \cite{Grubmuller:1994wj} had introduced the concept of conformational substates in 1994. \cite{Deuflhard2003,Schutte2003} provided some surveys on this topic.

Central to this model is the concept of transfer operator and the study of its eigenvectors and eigenvalues. Relevant reaction rates (as well as mean passage times, mean exit times, \ldots) can then be extracted from these eigenvalues. We provide a brief account of the key mathematical objects. The definitions can be made in a relatively general context (\cite{Schutte2003}). We consider a stochastic transition kernel $p(t,x,A)$ such that:
\begin{equation}
	p(t,x,A) = \mathbf{P}[X_{t+s} \subset A \, | \, X_s = x]
\end{equation}
where the family $\{X_t\}$ is a homogeneous Markov process indexed by a time variable $t$; $p(t,x,A)$ is therefore the probability that a Markov process started at $x$ is in $A$ after a time span $t$. This allows defining the Perron-Frobenius operator $P_t$ (propagator or forward transfer operator). See \cite{Schutte2003,Schutte:1999vh} for a definition that does not assume that the transition kernel $p(t,x,y)$ is continuous with respect to a probability measure $\mu(dx)$. Then:
\begin{equation}
	P_t \; u(y) = \frac{\displaystyle \int p(t,x,y) u(x) \rho(x) \, dx}{\rho(y)}
\end{equation}
The density $\rho$ is assumed to be invariant. In a similar fashion, ensemble transition probabilities can be defined as:
\begin{equation}
	p(s,C,D) = \frac{1}{\rho(C)} \int_C p(s,x,D) \rho(x) \, dx
\end{equation}

\subsection{Perron cluster cluster analysis}

The metastable sets alluded to earlier are defined as sets that are almost invariant under the Markov process. With our definition, the set $C$ is almost invariant if:
\begin{equation}
	p(s,C,C) \approx 1
\end{equation}

These almost invariant sets can also be identified by considering the eigenvalues and eigenvectors of $P_t$. For example, the density $\rho$, being invariant, satisfies $P_t \chi_\Omega = \chi_\Omega$ ,where $\Omega$ is the space of all conformations of the molecule and $\chi_\Omega$ is the characteristic function of $\Omega$ (in that case simply the function $\boldsymbol{1}$). The eigenvalue 1 therefore corresponds to the equilibrium distribution $\rho$.

Eigenvalues close to one form the so-called Perron cluster. The largest of these eigenvalues (not equal to 1) can then be associated with the slowest rate in the system. See the thesis of \cite{Huisinga2003} for an extensive discussion of transfer operators and metastability. \cite{huisinga2006metastability} provides some upper and lower bounds on these eigenvalues. There is a significant body of literature on the identification of these sets, in particular using an approach called Perron cluster cluster analysis (PCCA). The premise is that if one has uncoupled Markov chains (that is $P_t$ is a block diagonal matrix assuming discrete states) then, assuming $k$ separate aggregates or clusters, we will have $k$ eigenvectors with eigenvalue 1. Each eigenvector is constant over each cluster and changes sign (positive, negative or 0) between clusters (\cite{Deuflhard2000}). In PCCA, the sign has been used to identify these clusters or aggregates.

In a real application however, we are not dealing with uncoupled Markov chains but rather weakly coupled chains, resulting in a single eigenvector with eigenvalue 1 and a cluster of eigenvalues near 1, \cite{Schutte:vy}. The identification of the clusters based on the sign structure (\cite{Deuflhard2000}) is then more difficult as the sign change is more progressive with a smooth change of value across the transition region, and the determination of the sign is more difficult as the eigenvectors may assume very small values whose sign ($+1$, $-1$, $0$) is then difficult to determine.

One remedy to this is to recognize, as will be further discussed later on, that the sign of an eigenvector entry is not the right quantity to look at. For example the second eigenvector $\rho_2$ with eigenvalue $1-\varepsilon$ is typically nearly constant with value say $\rho_2^A>0$ in one cluster and is constant with value $\rho_2^B<0$ in another cluster. The transition region is not associated with the point where $\rho_2$ is zero (or changes sign) but rather where $\rho_2$ assumes the value $1/2 (\rho_2^A + \rho_2^B)$ (midpoint between the two plateaus). This midpoint value can be used in a robust fashion to determine the cluster boundaries.

Another approach pioneered by Deuflhard (\cite{deuflhard2005robust,Weber2004}) called robust Perron cluster cluster analysis (PCCA+) proposes a fuzzy decomposition where instead of a strict partitioning into clusters one calculates of partition of unity $\sum_i \tilde{\chi}_i(l) = 1$ (at a discrete state indexed by $l$) where each function $\tilde{\chi}_i$ is called an almost characteristic function which smoothly transitions from 1 to 0 outside out a cluster. In effect, \cite{Weber2004} assign a grade of membership between 0 and 1 to each state (in a discrete setting). Therefore, each state $l$ may correspond to different clusters with a different grade of membership, defined by $\tilde{\chi}_i(l)$. This approach was used to study a relatively long polyalanine (Ala$_8$ and Ala$_{12}$) in \cite{Noe2007}. The backbone torsion rotamer pattern was used to define the microstates.

The thesis of \cite{Weber2007} reviews these methods and discusses a meshless approach in which the membership functions are used to construct a Galerkin approximation of $P_t$. The fact that the basis functions are non-orthogonal (in contrast with partition functions for example) leads to a non-diagonal mass matrix in the Galerkin formulation.

We mention an alternative method based on a singular value decomposition, see \cite{Fritzsche2008}. The main drawbacks of the approaches mentioned above (PCCA and variants) are the difficulty of identifying the Perron cluster if the transition matrix of the Markov chain has no significant spectral gaps; in addition, the calculation of the eigenvectors may be badly conditioned if the Perron cluster contains many eigenvalues very close to 1. The SVD approach attempts to mitigates some of these issues.

Kube et al. (\cite{Kube:2005wv,Kube2007}) have used this decomposition into metastable sets with PCCA+ to construct a coarse grained matrix that approximates the exact fine grained propagator. Limitations of the resulting propagator are discussed.

An application of this approach to systems described by a Langevin equation is given in \cite{Schutte:2000tz}. An issue in this approach is the fact that the conformational space needs to be discretized appropriately in order to form a discretized approximation of $P_t$ (called transition matrix in \cite{Cordes2002}). This is difficult in practice since molecular systems live in high dimensional space. \cite{Cordes2002} propose to focus on dihedral angles to reduce the dimension and further improve their approach by considering a recursive decomposition in which space is first decomposed using the dihedral angle with the longest auto-correlation (which is shown to correlate with metastability). This leads to a first partitioning of space. Then, these metastable subspaces are further decomposed by applying the same strategy recursively (see \cite{Cordes2002}). PCCA is finally applied to the resulting coarse decomposition. \cite{galliat2002automatic} proposes an approach based on self-organized neural networks, also to attack this problem of dimensionality. \cite{schultheis2005extracting,Kloppenburg:1997gy} approach the problem of dimensionality through the use of density-oriented discretizations that represent the probability density using a mixture of normal distributions.

The concept of metastability also allows viewing the Markov chains as a collection of rapidly mixing chains within certain parts of the state space, that are weakly coupled together. This gives rise to the concept of uncoupling wherein uncoupled Markov chains, which resemble the original chain, are formulated for each metastable states. Then these $k$ chains can be coupled again by a $k \times k$ coupling matrix. The resulting system contains all the important information from the original chain. This is discussed in \cite{meerbach2005eigenvalue}, where in addition upper bounds are provided for the 2nd eigenvalue of the uncoupled Markov chains to establish that they are indeed rapidly mixing.

In a related approach to building coarse models without resorting to the transition matrix or transfer operator, \cite{Horenko2006} attempts a direct construction of a coarse grained model using a hidden Markov model with output given by stochastic differential equations (SDE). The model is then formulated as the combination of a set of SDEs (Langevin equations) and a rate matrix which determines how one jumps from an SDE to the next. See \cite{Horenko2006,Horenko2008a}. The idea of hidden Markov model is reused in \cite{fischer2007identification} but this time the output of a hidden Markov state is a probability density function in the observed variables (in this case torsion angles along the protein backbone). Von Mises output distributions are used since the observed variables are angles.

Although many clustering methods have been proposed, we mention the work of \cite{coifman2008diffusion} who proposed to build a coarse grained representation based on the eigenvectors of the diffusion map. The algorithm is based on the definition of a weighted graph on the simulated points and the subsequent computation of the first few eigenvalues and eigenvectors of a random walk on this graph. Connections are made to the backward Fokker-Planck operator.

A related approach was developed by Shalloway and his group. See \cite{Oresic:1994ij, Anonymous:1996wx, Shalloway:1996en, Ulitsky:1998bt, Korenblum:2003dl}. \cite{Church:1999fl} review these methods. In \cite{Oresic:1994ij}, Gaussian packets are used to characterize metastable basins and equations are provided to calculate and evolve packets. Packets are Gaussian functions in $\Omega$. These packets can be numerically obtained by computing the average value of $x$ locally in some metastable basin, and the variance of $x$:
\begin{equation}
	K_i^{-1} = 2 \beta \; \big\langle (x-x^0_i)(x-x^0_i)^T \big\rangle_i
\end{equation}
where $x^0_i$ is such that $\langle K_i \, (x-x^0_i)\rangle_i = 0$. See \cite{Oresic:1994ij} for the notation $\langle \; \rangle_i$ (this is a suitable local averaging in basin $i$). Then the Gaussian packet has the form:
\begin{equation}
	p_i(x) = \exp \big(-\beta [ V_i + (x-x^0_i)^T K_i (x-x^0_i)] \big)
\end{equation}
In \cite{Shalloway:1996en}, the eigenfunction expansion of the Smoluchowski equation is used to optimize the macrostate expansion (using Gaussian packets) by satisfying a minimum uncertainty condition. \cite{Ulitsky:1998bt} builds macrostates based on a variational principle, from which the transition region (separating the macrostates) can be identified and analyzed. Connections with the variational transition state theory (Pollak in \cite{fleming1993activated}, p. 5--41) are made.

\subsection{Markov state models}

\label{sec:msm}

Another broad class of methods are Markov state models (MSM). Markov chains have a long history, however their application to bio-molecular modeling and protein modeling is relatively recent and goes back to papers by \cite{Singhal:2004is}, \cite{Swope:2004gx}, and \cite{Swope:2004da}. Two landmark papers are \cite{chodera2007automatic} and \cite{chodera2006long}. See~\cite{pande2010everything} for a review and discussion of this model. This approach is closely related to conformation dynamics and many theoretical results from conformation dynamics directly apply to MSM. To build a discrete model of the continuous underlying conformational space $\Omega$, the method starts by partitioning $\Omega$ into cells. The advantage is that this decomposition is often easier to obtain in practice than defining a reaction coordinate. In this respect this approach has some of the advantages of transition path sampling, which also does not require a precise knowledge of the reaction coordinate. This is a significant departure from TIS, FFS, and milestoning who rely primarily on a single order parameter that measures the progress of the reaction from $A$ to $B$ in an essentially sequential way (progress along a single dimensional variable or order parameter). MSM remains easy to construct even in the presence of multiple pathways with no obvious reaction coordinate.

These macro-states can be constructed in different ways, for example from pathway sampling information (\cite{Singhal:2004is}) or by partitioning the Ramachandran map (\cite{chodera2006long}). \cite{Schutte1999} uses a decomposition based on torsion angles. In this paper, a method is also proposed to define generalized angle coordinates, in an attempt to reduce the number of coordinates in the problem. It is based on ideas from \cite{Amadei1993}, and statistical analysis of circular data (\cite{Fisher:1993uw,Fisher:1983kt}). \cite{chodera2007automatic} proposes an automatic procedure to create macro-states using an iterative procedure and the $k$-medoids algorithm (a partitioning algorithm similar to the $k$-means algorithm). Microstates are iteratively lumped into macrostates using the $k$-medoids algorithm, and split again into microstates to iteratively refine the definition of the macrostates. 

\begin{figure}
\centering
    \beginpgfgraphicnamed{tikzfigs/fig22-publisher}%
    \input{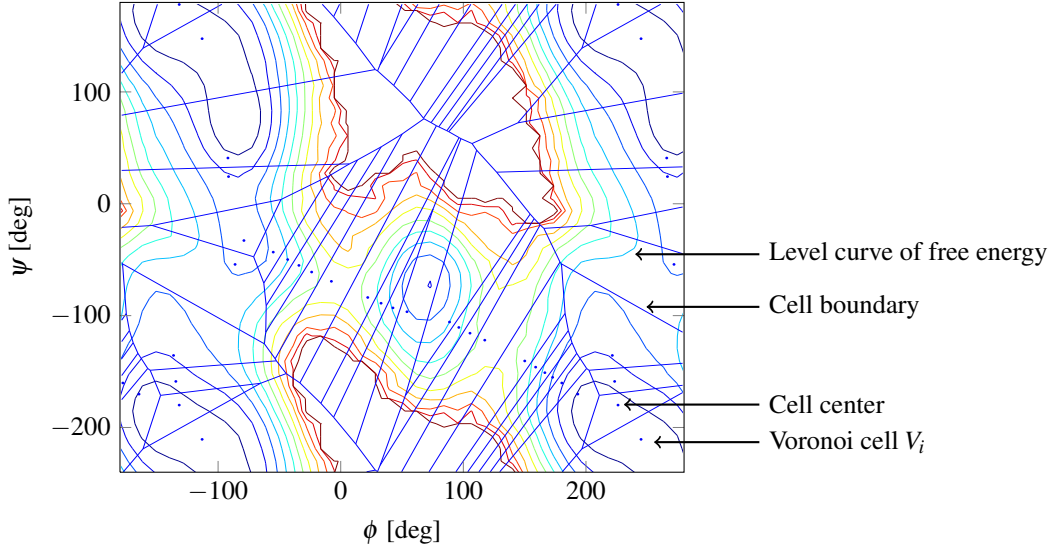}%
    \endpgfgraphicnamed%   
  
\caption{Example of Voronoi cells $V_i$ for alanine dipeptide. The angles $\phi$, $\psi$ are torsional angles along the backbone of alanine dipeptide. The cell centers are shown using small blue dots. The blue lines denote cell boundaries (milestones $S_i$). In this case, cell boundaries are by construction equidistant to two centers. The colored curves are level curves of the free energy. Low energy regions are dark blue while high energy regions are dark red. \label{vi}}
\end{figure}

Voronoi tesselation has also been proposed since it allows a simple construction of macro-states based only on the definition of the cell centers (\cite{VandenEijnden:2009dw}). See Fig.~\ref{vi}. Given a set of points $\{x_i\}$, a Voronoi cell $V_i$ is defined as:
\begin{equation}
	V_i = \{ x \; | \; |x-x_i| \le |x-x_j|, \quad j \neq i \}
\end{equation}
Such cells are convenient since there is a straightforward equation to determine in which cell a given point $x$ is. The centers $x_i$ can be obtained in different ways including simulations at high temperature, using nudge elastic bands or the string method, techniques to sample transition pathways, or computing minimum (free) energy pathways. We will show later on that placing centers along minimum energy pathways is often a good choice (see \cite{Pan:2008bp} for a related discussion).

We note that in \cite{schultheis2005extracting} this issue is altogether avoided by constructing a density-oriented discretization (an attempt to circumvent the curse of dimensionality) based on approximating the probability density in conformational space using a mixture of univariate normal distributions.

Once the macro-states have been defined, one calculates the transition matrix:
\begin{equation}
	P_{ij}(\tau) = \mathbf{P}(\text{particle in cell $i$ at time 0 is in cell $j$ at $\tau$})
\end{equation}
where $\tau$ is the so-called lag-time. The eigenvalues and eigenvectors of this matrix allow computing the different rates in the system, along with mean passage times (\cite{Swope:2004gx,chodera2006long}).

\cite{Park:2006io} use the concept of Shannon entropy to measure and identify non-Markovity. Shannon entropy measures the amount of uncertainty associated with a random variable. Non-markovity can be measured by evaluating the change in uncertainty (Shannon entropy) for a Markov variable $X_n$ if one prescribes the previous state $X_{n-1}$ vs.\ the last two states $X_{n-1}$ and $X_{n-2}$. Markov chains are such that the uncertainty is the same in both situations. From the definition of the Shannon entropy
\begin{align}	
	\text{If the previous state is known:} \quad &
	H(X_n|X_{n-1}) = - \sum_{x,y} \mathbf{P}(X_n=x,X_{n-1} = y) \; \ln \mathbf{P}(X_n=x \, | \, X_{n-1} = y) \\
	\text{If the last two states are known:} \quad &
	H(X_n|X_{n-1},X_{n-2}) = \\
	& \hspace*{-0.3in} - \sum_{x,y,z} \mathbf{P}(X_n=x,X_{n-1} = y,X_{n-2} = z) \; \ln \mathbf{P}(X_n=x \, | \, X_{n-1} = y, \, X_{n-2} = z) \\
	\text{\bf Measure of non-Markovity:} \quad & 
	R = \frac{H(X_n|X_{n-1})-H(X_n|X_{n-1},X_{n-2})}{H(X_n|X_{n-1})} 
\end{align}
From this definition, $R=0$ for Markov chains\footnote{$H(X_n|X_{n-1},X_{n-2}) \le H(X_n|X_{n-1})$ so that $R \ge 0$. Note that $R=0$ does not imply that the system is Markovian since it is possible that $R=0$ and $H(X_n|X_{n-1},X_{n-2},X_{n-3}) < H(X_n|X_{n-1})$.} and $R=1$ if the knowledge of $X_{n-2}$ and $X_{n-1}$ completely determines $X_n$. This procedure can also be used to refine the definition of macro-states. Another approach to identify and correct non-Markovity is proposed in \cite{Nerukh:2010ix} using the concepts of ``computational mechanics'' (a term coined by \cite{Crutchfield:1994gh}) and $\epsilon$-machines (\cite{Crutchfield:1989bx,Shalizi:2001in}).

\cite{Huang2009} discusses the use of generalized ensemble algorithms, e.g., the replica exchange method, parallel tempering or simulated tempering (\cite{Hansmann:1999va, Sugita:1999cl, Lyubartsev:1992js, Marinari:1992tl}) and how they can be combined with Markov State models to provide an efficient conformational sampling algorithm. See \cite{Bowman2009} for a discussion of similar ideas and how they have been implemented in the piece of software MSMBUILDER.

\cite{Chiang:2010bq} proposes to construct Markov models with hidden states as a way to construct more accurate models with fewer states (hidden Markov models). In such models, the hidden Markov states do not represent individual protein conformations but rather overlapping probabilistic distributions over the conformational space.

Applications of this approach are discussed in several papers including:
\begin{itemize}
	\item a polyphenylacetylene (pPA) 12-mer in explicit solvent for four common organic and aqueous solvents (acetonitrile, chloroform, methanol, and water): \cite{Elmer:2005gt,Elmer:2005gy}.
	\item lipid membrane fusion: \cite{Kasson:2006hl}.
	\item the villin headpiece: \cite{Jayachandran:2006de}. \cite{Bowman:2009jw} discusses the application of the software MSMBUILDER to the villin headpiece.
	\item polyalanines: \cite{Noe2007}. See also \cite{Noe:2008bg} with a review and discussion of MSM.
	\item PinWW domain: \cite{2009PNAS..10619011N,Morcos:2010hn}. In addition, in this paper, it is shown how folding pathways can be reconstructed from the MSM using transition-path theory (\cite{E:2006hg,Metzner2009}).
	\item the millisecond folder NTL9 (1--39): \cite{Voelz:2010hs}.
\end{itemize}

\section{Non-equilibrium umbrella sampling and reactive trajectory sampling}

\subsection{Non-equilibrium umbrella sampling}

In recent years, some approaches have in some sense tried to combine ideas from the previous sections, for example by calculating reactive trajectories or transition pathways from $A$ to $B$ (thereby being exact, contrary to Markov state models whose accuracy depends on the Markov assumption), while using a decomposition of the conformational space $\Omega$ into cells. These methods therefore combine the advantages of transition path sampling with the sampling efficiency of MSM, which requires only ``local'' sampling.

The first approach we will discuss is the one proposed by Dinner and co-workers, \cite{Warmflash:2007dz, Dickson:2009gt}. The method originates from the method of umbrella sampling, \cite{Torrie:1977wa}, in the sense that it tries to enhance sampling in poorly sampled region. The latter was broadly speaking adapted to allow modeling non-equilibrium systems. In this case, when computing a reaction rate from $A$ to $B$, particles are removed from the system each time they reach $B$ and are reinjected in $A$, thereby creating a steady-state but out of equilibrium situation.

This technique constructs two staggered lattices (using boxes, \cite{Warmflash:2007dz}, or following minimum energy pathways, \cite{Dickson:2009gt}). Simulations are run inside each macro-state (a box in \cite{Warmflash:2007dz}). Then one records when the system attempts to leave a box. At that point, two strategies are applied.

First, we keep track of the number of particles going from box $i$ to $j$. Each time a particle attempts to go from $i$ to $j$, the weight of box $i$, $W_i$, and $j$, $W_j$ are adjusted according to:
\begin{equation}
	- \Delta W_i = \Delta W_j = s \, W_i \, \frac{T^*}{T_i} 
	 \label{eq:weight}
\end{equation} 
where $T_i$ is the time elapsed in region $i$ (to account for situations in which longer simulations are run in some boxes), $T^*$ is some arbitrary time scale to make the equation dimensionally correct, and $s$ is a small parameter used to adjust the rate at which the weights $W_i$s vary. With this equation, the weight of each box converges to its correct steady-state value.

Second, one needs to determine which point should be used to reinsert a walker that left box $i$ back into box $i$. For this walkers that leave another box $j$ and attempt to enter $i$ are saved and, among those, one is picked according to the probability rule
\begin{gather}
	p_a = \frac{\displaystyle \sum_b N_{ba} W_{j(b)} / T_{j(b)}}{Z_i} \label{eq:reinsert}\\
	Z_i = \sum_{ab} N_{ba} W_{j(b)} / T_{j(b)} \notag
\end{gather}
where $p_a$ is the probability of choosing state $a$ in box $i$, $N_{ba}$ is the number of crossings from state $b$ in box $j(b)\neq i$ to state $a$.

In \cite{Warmflash:2007dz}, page 154112-4, end of section G, it is argued that two lattices are needed. Simulations are run in both but walkers that re-enter box $i$ in lattice 1 are chosen with Eq.~\eqref{eq:reinsert} using data from lattice 2. Otherwise it is claimed that the method is unstable and convergence may not be achieved. The argument put forward is as follows:
\begin{quote}
Suppose, for example, that the weight of a box ($B$) fluctuates upward. By Eqs. (3) and (4) [Eq.~\eqref{eq:reinsert} in this manuscript], walkers in neighboring boxes will then be reset to boundary states accessible from $B$ more often. However, if transitions from those states to ones in $B$ are allowed, with some probability, the reset walkers will immediately attempt to enter $B$ and increase its weight further according to Eq. (5) [Eq.~\eqref{eq:weight} in this manuscript]. This positive feedback loop causes the single-lattice scheme to be unstable in simulations to obtain the steady-state probability distribution as a function of multiple variables. The use of two lattices enables boundary states on one lattice to be chosen using the fluxes from the other lattice, which breaks the feedback loop and enables convergence.
\end{quote}
This argument is not so clear unfortunately. As a particle leaves $i$ to enter $j$, the weight $W_i$ is reduced by $-s W_i T^*/T_i$. If the particle reenters $i$, the weight is increased again by $s W_j T^*/T_j$. In general if $W_i$ has an upward fluctuation, the net result is a reduction in $W_i$. In some cases, walkers do not re-enter $i$ and move to some other box, further reducing $W_i$. As a result, in the scenario mentioned above of a temporary fluctuation upward of the weight, on average, Eq.~\eqref{eq:weight} will slowly reduce the weight of box $i$ and return it to its correct steady-state value.

Following \cite{Warmflash:2007dz} and \cite{Dickson:2009gt}, \cite{VandenEijnden:2009jh} developed a similar approach but that uses a single partition of space, based on Voronoi cells. The approach is similar with the following differences:
\begin{itemize}
	\item Weights are adjusted based on fluxes between cells. These fluxes, $N_{ij}/T_i$ (where $N_{ij}$ is the number of crossings from $i$ to $j$), are used to solve a linear system that provides an approximation to the steady-state weights $W_i$:
	\begin{equation}
		\sum_{j,i \neq j} W_i \; \frac{N_{ij}}{T_i}  = \sum_{j,i \neq j} W_j \; \frac{N_{ji}}{T_j}
	\end{equation}
	Another ``global'' scheme to adjust the weights (by contrast with the local scheme~\eqref{eq:weight}) is given in \cite{Dickson:2009fu}.
	\item To pick a re-entry point, one of the {\it boundaries} $\alpha$ of cell $i$ is randomly picked using a probability law obtained from the cell flux and steady-state probabilities:
\begin{equation}
	p_{\text{boundary $\alpha$ of cell $i$}} = \frac{W_j \; N_{ji} / T_j}{
	\sum_{k,k\neq i} W_k \; N_{ki} / T_k
	} 
\end{equation}
where boundary $\alpha$ is the boundary between cell $i$ and $j$.	
Although the implementation is different, this is similar to \cite{Warmflash:2007dz}.
\end{itemize}

Based on this, it appears that using two lattices is not necessary and that the scheme correctly works with a single lattice. Also the idea of lattice is no longer discussed in a more recent publication, \cite{Dickson:2010gf}. See also \cite{Dickson:2011vl} for an application of this method to unfolding and refolding of RNA. In this paper as well, a single lattice is used. \cite{Dickson:2010gf} present some theoretical results regarding non-equilibrium umbrella sampling (an analysis of the convergence of the weights using the local scheme), a comparison with and discussion of forward flux sampling, and recent applications of these methods.

\subsection{Reactive trajectory sampling}

The second method, which is related in some fashion to the previous class of techniques, can be attributed to \cite{Huber:1996dn}. In this reference, the method is developed assuming that an approximate reaction coordinate has been chosen. However, it is not difficult to extend this approach to a general decomposition of the conformational space $\Omega$ in a manner similar to, for example, \cite{VandenEijnden:2009dw} with Voronoi cells. This method will be discussed in more details below. It consists in running a large number of simulations (or ``walkers'') in parallel in such a way that a given number of walkers are maintained in each cell or macro-state. Macro-states that are near an energy barrier will tend to be depleted and therefore a strategy is applied to duplicate walkers in this macro-state, in a statistically correct way. This is done by assigning statistical weights to each walker. For example a walker with weight $w$ can be split into two walkers, starting from the same location in $\Omega$, with weights $w/2$. Conversely, macro-states that are at low energy will tend to become overcrowded and walkers are then removed. If for example we have two walkers with weights $w_1$ and $w_2$, we randomly select one with probabilities $(w_1/(w_1+w_2),w_2/(w_1+w_2))$ and assign to it the weight $w_1+w_2$. This approach ensures an efficient sampling of phase space.

In order to calculate a reaction rate, the macro-state corresponding to region $B$ is transformed into a cemetery state, that is any walker that enters this macro-state is removed from the simulation and re-inserted in region $A$. In this fashion, although the simulation is effectively out of equilibrium, the population of walkers is kept constant. This method allows computing all the relevant quantities of interest, such as reaction rates, free energy, metastable states, etc. We note that contrary to Markov state models, this approach does not suffer from non-Markovity errors and that in the limit of infinite sampling it provides an exact answer.

In \cite{Zhang:2007jc}, this technique was applied to explore the transition paths ensemble in a united-residue model of calmodulin. See also \cite{BinWZhang:2009tr, Bhatt:2010ea}. In \cite{Zhang:2010kf}, it is shown that the method initially developed in \cite{Huber:1996dn} is really applicable to a much wider class of problems and proposes some generalizations of this procedure.

We mention that a similar technique has been applied to simulated annealing to find minima of rough (or even fractal) functions (see \cite{Huber:1997ho}).

\bigskip

{\bf Detailed discussion of reactive trajectory sampling.} We now discuss in more details the method of \cite{Huber:1996dn,Zhang:2007jc,Zhang:2010kf} which we rename reactive trajectory sampling method (RTS), in the broader context of macro-state models (e.g., Voronoi cell partitioning). In this approach, systematic errors arising from non-Markovian effects are avoided by directly calculating reactive trajectories from $A$ to $B$ and obtaining the probability flux entering $B$ (or $A$ for the backward rate), \cite{Metzner:2006du}. When the energy barrier is high, this can be very inefficient since very few trajectories (if any) will make it to $B$ when started from $A$. However a simple trick allows improving the efficiency of the calculation to the extent that the decay of the statistical errors becomes essentially independent of the energy barrier height.

As before we split the space of possible configurations into cells. Then a large number of random ``walkers'' are initialized and advanced forward in time. The basic idea is to use a strategy whereby, in cells that get overcrowded (too many walkers), we merge walkers, thereby reducing their numbers, while in cells that are depleted (near transition regions), we split walkers to increase their number. The end goal is to maintain a given target number of walkers in each cell. With such an approach we are able to observe a constant stream of walkers going from $A$ to $B$ (and vice versa) irrespective of the height of the energy barrier. We now explain the details of the method.         

Assume we have $\nw$ walkers whose position gets updated at each time step. It is possible to resample from these walkers without introducing any bias in the calculation using the following procedure. Each walker, whose position is denoted $x_i$, is assigned a probabilistic weight $w_i$, for example initially equal to $1/\nw$. A walker can be split into $p$ walkers with weight $w_i/p$. After the split, each walker can be advanced independently. Averages can then be computed using:
\begin{equation}
	\langle f \rangle = \lim_{\nw \to \infty} 
	\frac{1}{\sum_j w_j}
	\sum_i w_i \; f(x_i)
	\label{eq15}
\end{equation}
This equation is always true irrespective of how many times the splitting procedure is applied, or how many steps are performed, as long as the initial position of the walkers is drawn from the equilibrium distribution. This is proved from the fact that the equilibrium distribution is by definition invariant under the dynamics under consideration for $x(t)$.

The reverse operation is possible. Assume we have a group of walkers with weights $w_1$, \ldots, $w_\nw$. Suppose we randomly pick a walker with probabilities $w_1/\sum_k w_k$, \ldots, $w_\nw/\sum_k w_k$, and assign to it a weight of $\sum_i w_i$. Since the average weight of walker $j$ is $(\sum_i w_i) \; (w_j/\sum_k w_k)$, Eq.~\eqref{eq15} remains true. This procedure can be used to reduce the number of walkers in a cell.

The algorithm below, called {\tt resample}, explains the procedure for resampling walkers. Before the procedure \verb|resample| is called we have several walkers in each cell with varying weights. The procedure \verb|resample| loops over the cells and select walkers in a way such that the new set of walkers all have the same weight, equal to the average weight of walkers in the cell. This is an important aspect of the method as assigning a constant weight can be proved to be optimal in terms of minimizing the variance and therefore statistical errors. This approach is different from \cite{Huber:1996dn}, which leads to walkers with varying weights and results in a somewhat larger variance. The proof is given below.

% \IncMargin{1em}
% \begin{algorithm}
% \caption{{\tt resample}. The algorithm is written in the programming language python. It was written by Eric Darve and Jes\'us Izaguirre (Notre-Dame University). \label{algo1}}

\lstset{language=python}
\lstset{numbers=left, numberstyle=\tiny, stepnumber=2, numbersep=5pt}
\lstset{keywordstyle=\color{black}\bfseries}
\begin{lstlisting}[frame=lines]
# This algorithm, called resample, is written in the programming language
# python. It was written by Eric Darve and Jesus A. Izaguirre (University
# of Notre-Dame).
# Input: list of walkers (list0) and list of walker weights (weights).
# Input: target number of walkers (ntargetwalkers).
# Output: list of walkers (list1) and their weights (newweights).
# weights[x] must be the weight of the walker with ID x.

from numpy import floor, argsort, random, sum
list1 = []       # new list of walkers
newweights=[]    # weights of new walkers
nwalkerlist1 = 0 # number of walkers in list 1

wi = # Initialize the list of weights for walkers in the current cell.
ind = argsort(-wi) 
# Sort the walkers in descending order based on their weights.
list0 = list(list0[ind])

W = sum(wi)
tw =  W / ntargetwalkers
# ntargetwalkers is the target number of walkers in cell

x = list0.pop() # We assume that there is at least one walker in the cell

while True: # while loop exits using a break.
    Wx = weights[x]
    if (Wx >= tw or len(list0) == 0):
        r = max(1, int(floor( Wx / tw )))
        # max is required because of round-off errors
        r = min(r,ntargetwalkers-nwalkerlist1) 
        # required because of round-off errors
        nwalkerlist1 += r         # update the number of walkers in list1		
        for item in repeat(x,r):  # insert r copies of walkers in list1
            list1.append(item)
            newweights.append(tw)
        if nwalkerlist1 < ntargetwalkers and Wx - r*tw > 0.0:
            list0.append(x)
            weights[x] = Wx - r*tw
        if len(list0)>0:
            x = list0.pop()
        else:
            break
    else:
        y = list0.pop()
        Wy = weights[y]
        Wxy = Wx + Wy
        p = random.random() # randomly select a walker
        if p < Wy / Wxy:
            x = y
        weights[x] = Wxy
\end{lstlisting}
   
% \DontPrintSemicolon
% \Output{An array, list1, of lists of walkers with constant weights in each cell}
% \BlankLine
% \For{i=0 \KwTo number of cells $-$ 1}{
% list0 = set of walkers in $V_i$\;        
% W = total weight of walkers in list0\; 
% tn = target number of walkers in $V_i$ (integer) \;
% tw = W/tn\;
% Sort list0 by weight in descending order\;
% list1[i] $= \emptyset$\;
% list2 $= \emptyset$\;
% Move a walker from list0 to list2\;
% \While{list0 $\neq \emptyset$}{
% W = total weight of walkers in list2\;
% \If{W $\ge$ tw}{
% r = $\lfloor$W/tw$\rfloor$\;
% Insert r copies of x in list1[i] with weight tw\;
% Insert one copy of x in list0 with weight W$-$r*tw\;
% list2 $= \emptyset$\;
% \If{list0 $\neq \emptyset$}{
% Move a walker from list0 to list2\;
% }
% }
% \Else{                        
% Move a walker from list0 to list2\;
% w = array of weights of walkers in list2\;
% Randomly pick a walker x in list2 with probabilities 
% \{w[0]/W, w[1]/W\}\;
% Replace the content of list2 by x with weight W\; 
% }}}
% \end{algorithm}

% \DecMargin{1em}  

We note that this algorithm terminates since when the last walker is removed from list0 we exactly have W = tw (target weight) so that the while loop does terminate. The maximum number of iterations in the while loop is bounded by $\nw\ + $ tn: the number of walkers in the cell before the procedure starts + the target number of walkers. (Lines 26 through 40 are executed at most tn times, while the lines 42 to 48 are executed at most $\nw$ times.) In addition it is apparent from line 33 that the walkers have the same weight at the end. The sorting of the weights on line 13 helps reduce data correlation. Indeed when a walker is split, samples become correlated for some time. The initial sorting makes sure that only walkers whose weight is greater than tw are split. As soon as we have processed all the walkers with weight greater than tw, r stays equal to 1. The reinsertion on line 35 is required to ensure a constant total weight W. The weight that is used, W$-$r*tw, ensures that the resampling is correct and that the total weight in each cell is unchanged by \verb|resample|.

\bigskip

\noindent {\bf Proof of optimality.} To simplify the discussion we assume that we have $\nw$ particles with weights $w_i$ such that $\sum_i w_i = 1$. We evolve the system in time such that the walker locations $\{x_i\}$ become uncorrelated. As a consequence, the weights $w_i$ are independent from the positions $\{x_i\}^\nw_{i=1}$. We also assume that the walkers do not have any particular order so that the statistics of $x_i$ and $w_i$ are the same as $x_j$ and $w_j$, $j \neq i$. The estimator of a particular quantity $f$ is
\[
\hat{f} = \sum^\nw_{i=1} w_i \; f(x_i)
\]
Then we have
\begin{equation}
\langle \hat{f} \rangle = 
\sum^\nw_{i=1} \langle w_i \rangle \; \langle f(x_i) \rangle
=
\langle \sum^{\nw}_{i=1} w_i\rangle \; \langle f \rangle
= \langle f \rangle
\end{equation}
which tells us that $\hat{f}$ is indeed an unbiased estimator of $f$. The statistical error can be estimated from the variance of $\hat{f}$:
\begin{align}
\langle (\hat{f}-\mu)^2\rangle
&=
\Big\langle \Big(
\sum^{\nw}_{i=1}w_i (f(x_i)-\mu)
\Big)^2 \Big\rangle\\
&=
\sum^{\nw}_{i=1} \sum^{\nw}_{j=1}
\langle w_i w_{j}(f(x_i)-\mu) (f(x_j)-\mu)
\rangle\\
&=
\sum^{\nw}_{i=1}
\langle w_i^2 \rangle
\langle (f(x_i)-\mu)^2 \rangle\\
& = \nw \; \langle w^2 \rangle \;
\langle (f-\mu)^2 \rangle
\end{align}
Therefore
\begin{equation}
	\langle (\hat{f}-\mu)^2\rangle = \langle (f-\mu)^2 \rangle
	\; \Big( \frac{1}{\nw} + \nw \; \text{Var}(w) \Big)
\end{equation}
since $\langle w \rangle^2 = 1/n_\text{w}^2$. In our algorithm the weights are kept constant (if tn in Algorithm {\tt resample} is constant) so that $\text{Var}(w) = 0$. The statistical error is therefore minimized. In \cite{Huber:1996dn}, the weights are not constant resulting in a larger statistical error.

$\square$

\bigskip

% As will be seen later on in the numerical results section, RTS is easier to use in practice. At small times it may suffer from systematic errors depending on the choice of the initial distribution of the particles. If the initial walkers are distributed randomly, the statistics will be biased until the walkers relax to their equilibrium distribution. Techniques to calculate a good initial distribution of walkers, or to correct the weights of walkers while the calculation is running to accelerate convergence are given in Section~\ref{section:coloring}.  

An important difference with the Markov model based on cells (coarse states) is that MSM must be run by construction with a known lag time $\tau$. Post-processing is then required to determine whether $\tau$ is large enough. After this, the simulation may need to be run again with a larger $\tau$ if it is found that memory effects are important. On the contrary, the convergence of RTS is easier to monitor. One simply needs to record the particles that reach $B$ and calculate the average flux. If the error is found to be too large, one can simply continue the simulation to accumulate more statistics, without losing the data already gathered.

In \cite{Huber:1996dn}, walkers that enter $B$ are re-inserted in $A$ thereby ensuring a steady-state system. It is possible to use a slightly different procedure where colors are given to walkers such that when a walker last entered $A$ its color is blue, while it is red if it last entered $B$. In effect the color changes from blue to red the first time the walker enters $B$ (similarly with $A$). The population for both colors is kept constant in each macro-state according to the algorithm \verb|resample|. This algorithm allows computing both the forward and backward rates, the free energy, and the equilibrium distribution of particles (by considering all particles, of any color).

\subsection{Optimal cells} 

\label{opt_cells}

Even though the method was shown to always converge to the correct answer, the rate of convergence, which depends on the rate at which particles transition from $A$ to $B$, depends on the choice of cells. We provide some guidelines to help make a good choice of cells. We note that, as explained previously, RTS is always unbiased and is exact in the absence of statistical errors. This is in contrast with milestoning which is exact only when the milestones are iso-surfaces of the committor function. In that case the milestones are called optimal since they minimize systematic errors. Here, RTS is unbiased. The optimal milestones in this context are the ones that minimize the statistical error. Even though the terminology is the same, the meaning is therefore quite different.

RTS ultimately amounts to sampling reactive trajectories from $A$ to $B$. Reactive trajectories are defined as trajectories that leave $A$ and reach $B$ without reentering $A$ at any point. These trajectories cluster around the minimum (free) energy paths, \cite{2006JChPh.125b4106M}. These paths correspond to reactive trajectories going from $A$ to $B$ with maximum likelihood, that is the probability density associated with this trajectory is maximum.

To discuss properties of minimum free energy pathways (MFEP), it is convenient to use generalized coordinates $(\xi_1$, \ldots, $\xi_p)$ to describe the system and use the free energy $A(\xi_1,\ldots,\xi_p)$. This is practically a more useful description and it removes degeneracies such as translation and rotation invariance. For example these generalized coordinates can be chosen as a set of internal coordinates describing the shape or structure of a molecule. If some information is available about the system we can reduce the number of such variables to focus on the variables of interest for the reaction at hand.

The generalized coordinates can have any units, for example \AA\ or deg. This indicates that some kind of non-dimensionalization procedure is required to work with $\xi$. This non-dimensionalization can be derived in different ways. We shortly describe how this can be done.

If one assumes for example a Brownian model for these variables:
\begin{equation}
	\frac{d\xi}{dt} = \nabla D - \beta \, D \, \nabla A(\xi) + R \; \eta(t)
	\label{eq16}
\end{equation}
where $R$ is such that $R \, R^T = 2 \, D$, $\beta = (kT)^{-1}$, and $\eta(t)$ is a random term with a normal distribution and variance 1, and $A$ is the free energy. The tensor $D$ is the diffusion tensor. Under simplifying assumptions this tensor can be approximated by:
\begin{gather}
	D \approx \beta^{-1} \tau_D M_\xi^{-1} \\
	[M_\xi^{-1}]_{ij} = \sum_k \frac{1}{m_k} \frac{\p \xi_i}{\p x_k}
	\frac{\p \xi_j}{\p x_k} 
\end{gather}
where $m_k$ is the mass of atom $k$; $\tau_D$ is a time scale associated with the rate of decay of the auto-correlation function for $d\xi/dt$. The tensor $M_\xi$ is non constant. However we will assume that it can be approximated by its statistical average and that its fluctuations can be neglected.

From Eq.~\eqref{eq16}, we can conclude that along the MFEP we must have that the tangent $d\xi/ds$ (where $s$ is some parameterization such as the arc length) is parallel to $D \nabla A$ or equivalently:
\begin{equation}
	\frac{d\xi}{ds} \propto M_\xi^{-1} \; \nabla A
\end{equation}
This result is somewhat counter-intuitive since we would expect $d\xi/ds \propto \nabla A$ but is a result of the metric associated with $\xi$ and defined by $M_\xi$. This suggests normalizing $\xi$ using:
\begin{equation}
	\tilde{\xi} = M_\xi^{1/2} \; \xi
\end{equation}
With the variables $\tilde{\xi}$ we have the expected relation, along the MFEP:
\begin{equation}
	\frac{d\tilde{\xi}}{ds} \propto \tilde{\nabla} A
\end{equation}
where $\tilde{\nabla}$ involves derivatives with respect to $\tilde{\xi}$.

Similarly when defining Voronoi cells, the correct distance to use should be consistent with the rate of diffusion and therefore the following distance must be used:
\begin{equation}
	||\Delta \xi||_{M_\xi} = \big( \Delta \xi^T \; M_\xi \; \Delta \xi)^{1/2} = ||\Delta \tilde{\xi}||_2
\end{equation}

We will now assume that we are using $\tilde{\xi}$ instead of $\xi$, but we will keep the notation $\xi$ for simplicity. Returning to the issue of reactive trajectories and optimal cells, we note that in each cell, walkers tend to accumulate in low energy regions. Consequently we can expect the method to be efficient (practically the statistical errors are small) whenever the regions around the MFEPs are well sampled, that is the low energy region in a cell should overlap as much as possible with the MFEP.

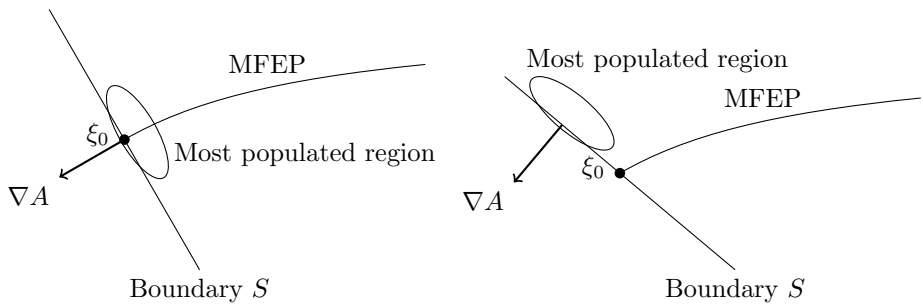
\begin{figure}
\centering	
\begin{tikzpicture}
\draw (0,0) .. controls +(30:1)	and (2,0.8) .. (4,1);
\fill (0,0) circle (2pt);
\path (0,0) node [anchor=east,xshift=-2pt,yshift=2pt] {$\xi_0$};
\path (1.9,1) node {MFEP};
\draw (-60:2) node [anchor=north] {Boundary $S$} -- (120:2);
\draw [->, thick] (0,0) -- (210:1) node [anchor=north east] {$\nabla A$};
\draw [rotate=-60] (90:0.2) ellipse (0.7 and 0.25);
\path (2.4,-0.2) node {Most populated region};
\end{tikzpicture} 
\begin{tikzpicture}
\draw (0,0) .. controls +(30:1)	and (2,0.8) .. (4,1);
\fill (0,0) circle (2pt);
\path (0,0) node [anchor=east,xshift=-2pt,yshift=2pt] {$\xi_0$};
\path (1.9,1) node {MFEP};
\draw (-40:2) node [anchor=north] {Boundary $S$} -- (140:2);
%\draw [->, thick] (0,0) -- (210:1);
\draw [rotate=-40] (-1,0.2) ellipse (0.7 and 0.25);
\path (0.5,1.5) node {Most populated region};
\draw [->, thick] (140:1) -- +(230:1) node [anchor=north east] {$\nabla A$};
\end{tikzpicture}
\caption{Schematic representation of a cell crossed by an MFEP. On the left panel, the boundary of the cell is normal to the MFEP, ensuring a good sampling along the MFEP. This is illustrated by the most populated region containing $\xi_0$. On the right panel, the orientation of the boundary was changed so that the most populated region is now shifted and no longer overlaps with the MFEP. In that case statistical errors are expected to be somewhat larger. In the extreme case where the boundary is parallel with $\nabla A$ the simulation converges very slowly.}
\label{fig:mfep}
\end{figure}

Consider a cell and assume that $A$ is not singular in that cell. This implies that $A$ is minimum on one of the boundaries, which we denote $S$. Let us assume that the MFEP crosses $S$ at $\xi^0$. In order to minimize the statistical errors, we impose the condition that $\xi^0$ is the point with the highest probability density in the cell. In that case, the boundary must be orthogonal to $\nabla A$ at that point, that is the tangent to the MFEP should be orthogonal to the cell boundary. As a note, we point out that along the MFEP, the gradient of the committor function is also parallel to $\nabla A$, which implies that locally the boundary $S$ is an iso-surface of the committor function. See Fig.~\ref{fig:mfep}.

As pointed out previously the method always converges to the correct answer. However the statistical errors are expected to increase as the cells degrade, that is the cell boundaries are no longer orthogonal to $\nabla A$. Then the walkers in a cell start accumulating in regions that are far away from the MFEP. As a result the probability of seeing a walker reach $B$ becomes correspondingly smaller, leading to larger statistical errors. Instead of observing a steady (with small fluctuations) stream of particles with small weights reaching $B$, we see a more sporadic stream of particles with larger weights. See Fig.~\ref{fig:mfep}.

\subsection{Metastability, convergence, and the multi-colored algorithm}

\label{section:coloring}

In RTS, walkers initially start from region $A$ and are moved out of the simulation when they reach region $B$. If nothing is done, the total weight slowly diminishes. When the reaction rate is slow, the simulation remains accurate. However, when multiple rates are comparable, this may lead to biased results. One option to prevent the total weight from diminishing is to reinsert a walker in basin $A$ after it enters basin $B$. Another possible strategy, as discussed earlier, with similar efficiency but somewhat easier to implement, is to use walkers with two different colors, for example blue and red. By definition, red (resp.\ blue) walkers are those that have last visited $A$ (resp.\ $B$). Whenever a blue walker enters $B$, its color changes to blue, and vice versa. The resampling algorithm is applied to each color separately so that each bin contains the same amount of walkers of each color. This approach allows computing the forward and backward rates, and the free energy.

The efficiency of RTS degrades when there are other slow reactions rates (values of $\lambda_k$ that are small) in the system, that is other metastable regions in $\bar{A}$ or $\bar{B}$, \cite{Noe2007}. In that case, the convergence of the algorithm is limited by the rate at which walkers equilibrate in region $\bar{A}$ and $\bar{B}$, which is given by $\lambda_3^{-1}$. This leads to long correlation times for the measured fluxes, and therefore overall slow decay of the statistical errors.

This can be addressed by assigning appropriate weights for the walkers at $t=0$. One possibility is to start running RTS for a while and calculate fluxes between cells, given by the weights of walkers moving from cell $i$ to $j$ after a time step divided by the total weight of walkers in cell $i$. From the flux matrix, the steady-state weights of cells can be computed. These weights can then be used to adjust the weights of walkers in each cell. In principle this needs to be iterated until convergence. However only one or two iterations are typically needed. The remaining deviations from steady-state can be later on reduced by the production RTS run. The goal of this procedure is simply to improve the initial weights in each cell in order to bypass the initial slow convergence.

The key property of this procedure is that the accuracy with which fluxes can be computed is independent of the presence of metastability and depends primarily on the sampling inside each cell. Overall convergence is a global property, which is why a few iterations may be required, but since this is followed by the production RTS run, this part of the calculation only requires low accuracy and few iterations. This approach is similar to a method described in \cite{Bhatt:2010df}, called ``Enhanced weighted ensemble attainment of steady state.''

Another approach can be applied to address this shortcoming. It consists in using more than $2$ colors. To explain this in a simple fashion we return to the milestoning framework of Section~\ref{milestoning}.
RTS can be reinterpreted in terms of the milestoning framework. The difficulty in the optimal milestoning approach described in Section~\ref{milestoning} is that the cells must be such that their boundaries are iso-surfaces of the committor function. Even though this is in general difficult to realize, there is a case for which the problem is simplified. Consider the following three cells: cell $V_0$ enclosing $A$, cell $V_1$ enclosing $B$.
%, and cell $V_2$, which is the complement: $V_2 = (V_0\cup V_1)^\text{c}$. We define $S_0$ (resp.\ $S_1$) as the interface between $V_0$ (resp.\ $V_1$) and $V_2$. 
% By definition of the committor function, these are exact iso-surfaces. Since we only have two surfaces, Eq.~\eqref{mpt} then reduces to a scalar equation, and it states that the mean passage time is equal to the inverse of the flux, which can be computed using RTS.
Let us assume that we have another metastable basin $C$. The same reasoning can be extended to an arbitrary number of basins. We define $V_2$ as a cell enclosing $C$ and $V_3$ as the complement: $V_3 = (\cup_{i=0}^2 V_i)^\text{c}$. We define $S_0$ as the interface between $V_0$ and $V_3$, and similarly for $S_1$, and $S_2$. In the multi-coloring algorithm, each cell is assigned a color, say 0 is blue, 1 is red, and 2 is green. Each time a blue particle enters cell $V_2$, its color changes to green, and similarly for the other colors. The total weight of all walkers is therefore constant. We use the resampling algorithm to maintain a constant number of walkers in each cell, separately for each color. We then measure the mean flux of particles turning blue to green, etc.

Eq.~\eqref{mpt} still holds for this approach. The fluxes $F_{ij}$, from cell $i$ to $j$, are simply related to the probability matrix $P$ through: $F = \Delta t^{-1} P$. Eq.~\eqref{eq14} is not needed for this approach. We can directly obtain $P$ from the flux values $F$, computed using RTS.

In general the milestone $S_2$ is not an iso-surface of the committor. However, following the proof that optimal milestoning gives exact rates, we will have proved that the rate with multi-coloring is exact if we show that: $p_{ij}$, the probability to cross $j$ after $i$, and $\tau_i$, the mean time before crossing another milestone, are independent of the previous milestone $k$ that was crossed (see page~\pageref{MSTproof} for the proof in the optimal milestoning case). This property is in fact true for the multi-coloring approach because of our choice for $S_2$ and the fact that the equilibration time for $V_2$ is very small. The key assumption is that $V_2$ must be associated with a metastable state (minimum energy basin) so that the relaxation time in $V_2$ is small compared to the mean escape time.

This shows that the mean passage time predicted using Eq.~\eqref{mpt} in the multi-coloring framework is for all practical purposes exact, with no systematic bias. This approach allows considering the case of multiple slow rates, with no significant degradation in efficiency. The computational cost merely grows with the number of colors, but is independent of the degree of metastability (the values of the first $\lambda_k$, $k=1$, 2, etc).

This approach has some conceptual similarities with the technique of core sets of \cite{2011JChPh.134t4105S}. Their analysis of accuracy (section E, pp. 204105-7) carries over to RTS with multiple colors. In addition, if ones applies the Galerkin discretization approach from \cite{2011JChPh.134t4105S} to RTS with $m+1$ colors, one can calculate the slowest $m$ rates in the system (or phenomenological rates, see \cite{2011JChPh.134t4105S} pp.~204105-4, section F) which may be of interest for certain applications. 

As a final note, we point out that the method is embarrassingly parallel and can be easily implemented on a parallel machine. This is important as this allows making only few changes to a serial (sequential) molecular dynamics code to make it run efficiently on a parallel cluster, with RTS, without having to parallelize the core of the code.

\section{Analysis of statistical errors in Markov state models}

We now discuss some mathematical results for Markov state models. In particular we will analyze the sensitivity of the eigenvalues to perturbations in the transition matrix. This analysis will lead to estimates for the statistical errors in the method. This will also lead to an analysis of the systematic errors, due to the finite lag time $\tau$ (the length of the short trajectories used to build the Markov state model).

Several papers have discussed error analysis in the context of conformation dynamics (\cite{sarich2010approximation,Prinz:2011id}) or Markov state models (\cite{Hinrichs:2007hj,Singhal:2005dk,Hinrichs:2007th}). \cite{2009PhRvE..80b1106M} takes a different approach to error analysis by considering a method that generates random transition matrices and as a result can estimate errors in various quantities computed from the transition matrix. In \cite{2009PhRvE..80b1106M}, it is argued that this approach is more accurate since it does not rely on Taylor expansions to approximate the impact of small variations of the numerical transition matrix from the exact matrix. In that sense this is a more direct estimate of the statistical errors and how they impact various quantities of interest (stationary distribution, eigenvalues, committor function, etc). \cite{2009PNAS..10610884D} discusses the effect of memory in building coarse grained models in the context of the Mori-Zwanzig formalism.

Although many results presented in this paper can be extended to more general stochastic equations, we are going to focus on the relatively simpler case of Brownian dynamics (\cite{1978JChPh..69.1352E}, p.~1355):
\begin{equation}
	dx(t) = \big( \nabla D_\text{B} - \beta D_\text{B}(x) \nabla U(x) \big) \, dt + R_\text{B}(x) \, dW(t)
	\label{eq1}
\end{equation}
where $W(t)$ is a Wiener process (see p.~66 in~\cite{Gardiner:1997tb}), $D_\text{B}$ is the diffusion tensor, $\beta^{-1} = k_\text{B}T$, and $R_\text{B}$ satisfies $R_\text{B}(x) R_\text{B}(x)^T = 2D_\text{B}(x)$.

\subsection{Eigenvectors and eigenvalues of the transition matrix}

The rate is typically computed by considering the eigenvalues of the transition matrix:
\begin{equation}
	P_{ij}(\tau) = \mathbf{P}(\text{particle in cell $i$ at time $0$ is in cell $j$ at time $\tau$})
\end{equation}
where we assume that the dynamics is given by~\eqref{eq1}. This matrix is basically used to construct a Markov state model of the system.
                   
From the stochastic equation~\eqref{eq1}, we can define the conditional probability $\rho(x,t|x_0,0)$, which is the probability to be at $x$ at time $t$ if the system was at $x_0$ at time 0. This probability can be expanded in terms of the eigenfunctions $\rho_k(x)$ of the forward Fokker-Planck equation (for Eq.~\eqref{eq1}), and the eigenfunctions $\psi_k(x_0)$ of the backward Fokker-Planck equation (\cite{Gardiner:1997tb}, p.~165):
\begin{equation}
	\rho(x,t|x_0,0) = \sum_k \psi_k(x_0) \; \rho_k(x) \; e^{-\lambda_k t}
	\label{eq2}
\end{equation}                                                           
where $\lambda_k$ are real and positive eigenvalues (the two sets of eigenfunctions are associated with the same eigenvalues). See p.~32 in \cite{schutte1999conformational}, p.~174105-5 in~\cite{Prinz:2011id}, and p.~166 in~\cite{Gardiner:1997tb}. Since
\begin{equation}
	\int \psi_k(x) \; \rho_{k'}(x) \; dx = \delta_{k,k'}
\end{equation}
$e^{-\lambda_k t}$ is an eigenvalue of the kernel $\rho(x,t|x_0,0)$:
\begin{equation}
	\int \rho(x,t|x_0,0) \; \rho_k(x_0) \; dx_0 = e^{-\lambda_k t} \; \rho_k(x)
\end{equation}
We will denote:
\begin{align}
	\langle \psi_k \rangle_i & 
	=
	\frac{\int_{V_i} \psi_k(x_0) \rho(x_0) \; dx_0}{
	\int_{V_i} \rho(x_0) \; dx_0} \\
	% \langle \psi_k \rangle_{j,k'} & =
	% \frac{\int_{V_j} \psi_k(x_0) \rho_{k'}(x_0) \; dx_0}{
	% \int_{V_j} \rho_{k'}(x_0) \; dx_0} \\  
	\rho_{jk} & = \Big[ \int_{V_j} \rho_k(x_{0}) \; dx_{0} \Big] 
\end{align}  
where $\rho (x)$ is the equilibrium distribution of the system. 

We sort the $\lambda_k$ in increasing order. For most systems, there is a single eigenvalue $\exp(-\tau \lambda_1)$ equal to 1 ($\lambda_1=0$) and the corresponding eigenvector is the stationary distribution $\rho(x)$. We are interested in estimating $\lambda_2$ by computing the eigenvalues $\mu_2$ of the matrix $P_{ij}$, and using $\lambda_k \sim -\ln(\mu_k)/\tau$. In general, $\lambda_2$ and $-\ln(\mu_2)/\tau$ differ leading to inaccurate estimates. However under certain assumptions, which will be discussed, $-\ln(\mu_2)/\tau$ provides an accurate estimate. In these circumstances, the Markov assumption made in building the model becomes accurate.

The term $P_{ij}$ can then be written as:
\begin{align}
	P_{ij} (\tau) & = 
	\frac{\int_{x_0 \in V_i} \int_{x \in V_j} \rho(x,\tau|x_0,0) \rho(x_0) \; dx \; dx_0}{
	\int_{V_i} \rho(x_0) \; dx_0
	} \\
	& =
	\sum_k e^{-\lambda_k \tau} \; \rho_{jk} \; \langle \psi_k \rangle_i
\end{align}
 
% Consider now, for a fixed $k$, the vector with components $\rho_{jk}$. If
% \begin{equation}
% 	\sum_j \langle \psi_{k'} \rangle_j \; \rho_{jk} = \delta_{k,k'} 
% 	\label{eq3}
% \end{equation}                                         
% then $e^{-\lambda_k \tau}$ is an eigenvalue. Indeed:
% \begin{equation}
% 	\sum_i \rho_{ik} \, P_{ij}
% 	= \sum_{k'} e^{-\lambda_{k'} \tau} \; \rho_{j,k'} \; \sum_i \langle \psi_{k'} \rangle_i \; \rho_{ik}
% 	= e^{-\lambda_k \tau} \; \rho_{jk}
% \end{equation}
% 
% We will use a few additional results which we now derive. Taking $t=0$ in Eq.~\eqref{eq2}:
% \begin{equation}
% 	\delta(x-x_0) = \sum_{k'} \psi_{k'}(x_0) \; \rho_{k'}(x)
% \end{equation}
% Integrate this equation against $\rho_{k}(x_0)$:
% \begin{equation}  
% 	\begin{split}
% 		\rho_{k}(x) = \sum_{k'} \rho_{k'}(x) \;
% 	\int \psi_{k'}(x_0) \; \rho_k(x_0) \; dx_0 \\
% 	\Rightarrow \qquad \rho_{k}(x) = \sum_{k'} \rho_{k'}(x) \;
% \sum_{j}\langle\psi_{k'}\rangle_{jk} \; \rho_{jk}
% 	\end{split} 
% \end{equation} 
% Therefore, we must have:
% \begin{equation}
% 	\sum_j \langle\psi_{k'}\rangle_{jk} \; \rho_{jk} = \delta_{k,k'}
% 	\label{eq4}
% \end{equation}
% This equation is always satisfied whereas Eq.~\eqref{eq3} requires additional assumptions to be valid. 

% Whenever Eq.~\eqref{eq3} holds for $k=2$, we are guaranteed that the rate computed using the eigenvalues of $P$ is exact. We therefore now discuss different assumptions for Eq.~\eqref{eq3} to hold for $k=2$. 

% In the following discussion we will assume that $\lambda_2 \ll \lambda_k$ for $k>2$.  

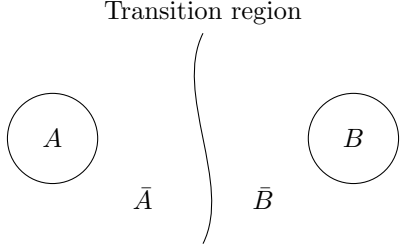
\begin{figure}
\centering
\begin{tikzpicture}[scale=2]
\draw (1,-0.2) .. controls (1.2,0.2) and (0.8,0.8) .. (1,1.2);
\path (0.6,0.1) node {$\bar{A}$};
\path (1.4,0.1) node {$\bar{B}$};
\draw (0,0.5) circle (0.3) node {$A$};
\draw (2,0.5) circle (0.3) node {$B$};
\path (1,1.2) node [anchor=south] {Transition region};
\end{tikzpicture}	
\caption{Definition of the different regions, $A$, $B$, $\bar{A}$, and $\bar{B}$.}
\label{fig:region}
\end{figure}

In the rest of this paper we sometimes have to make a distinction between the two minimum energy regions $A$ and $B$ and a partitioning of the space $\Omega$ into two metastable regions $\bar{A}$ and $\bar{B}$ (see Fig.~\ref{fig:region}). Typically $A$ is defined as a small region around a stable conformation of interest (the reactant state), and similar for $B$ (the product state). The regions $\bar{A}$ and $\bar{B}$ are defined as metastable regions, that is the rate of transition between these regions is the smallest among all other pairs of sets.

\bigskip

\noindent {\bf Properties of the eigenvectors.} 
The second eigenvector $\rho_2(x)$ can be used to define a partition  into two regions $\bar{A}$ and $\bar{B}$, with $A \subset \bar{A}$, $B \subset \bar{B}$. We will provide a more rigorous analysis later on but roughly speaking, the function $\rho_2/\rho$ is nearly constant in two regions, which define $\bar{A}$ and $\bar{B}$. This function changes sign between these two regions. The narrow region where $\rho_2/\rho \sim 0$ defines the transition region between these two metastable regions. See for example Chapter 6, p. 91--119 in~\cite{hill2004applied}. The molecule is assumed to have a high probability of being in region $A$ when in $\bar{A}$, and similarly for $B$.

\begin{figure}[htbp]
	\centering
 \begin{tikzpicture}

 \filldraw [fill=blue!20] (0,-0.2) rectangle (4,0.2);
 \filldraw [fill=green!20] (4,-0.2) rectangle (8,0.2); 

 \draw[->] (-1,0) -- (8.3,0);
 \draw[->] (0,-0.5) -- (0,2);
 \draw (8.6,0) node {$x$};

 \draw [thick] (0,1.5)  .. controls (7,1.5)  and (1,-0.4)  .. (8,-0.4);

 \draw (0,0)  .. controls (1,0)  and (1,1)  .. (2,1);
 \draw (2,1)  .. controls (3,1) and (3,0) .. (4,0);
 \draw (4,0)  .. controls (5,0) and (5,2) .. (6,2);
 \draw (6,2)  .. controls (7,2) and (7,0) .. (8,0);

 \draw (2.5,1.7) node {$\psi_2=\rho_2/\rho$};
 \draw (7.15,1.7) node {$\rho(x)$};

 \draw [dashed] (4,-1.4) -- (4,2);

 \draw [dotted,thick] (-0.3,1.5) node [anchor=east] {$(\rho(\bar{B})/\rho(\bar{A}))^{1/2}$} -- (-0.05,1.5);
 \draw [dotted,thick] (-0.3,0.55) node [anchor=east] {Midpoint} -- (4,0.55);
 \draw [dotted,thick] (-.3,-0.4) node [anchor=east] {$-(\rho(\bar{A})/\rho(\bar{B}))^{1/2}$} -- (8,-0.4);

 \draw [thick,->] (2,-1) node [anchor=north east,text width=2.8cm,xshift=1.4cm] {Transition point $\pi(x)=1/2$} -- (3.9,-0.06);
 \draw [thick,->] (5,-1) node [anchor=north,text width=4cm,xshift=1.5cm] {Points go to set $\bar{A}$ with probability $\rho(\bar{A})$}
	-- (4.53,-0.06);
	
 \path (2,2.4) node {Region $\bar{A}$};
 \path (6,2.4) node {Region $\bar{B}$};

 \filldraw (4,0) circle (0.05cm);
 \filldraw (4.5,0) circle (0.05cm);
 \filldraw (4,0.55) circle (0.05cm);

 \end{tikzpicture}
	\caption{Definition of the metastable regions $\bar{A}$ and $\bar{B}$ and the second eigenvector $\rho_2$.}
	\label{fig:eig}
\end{figure}
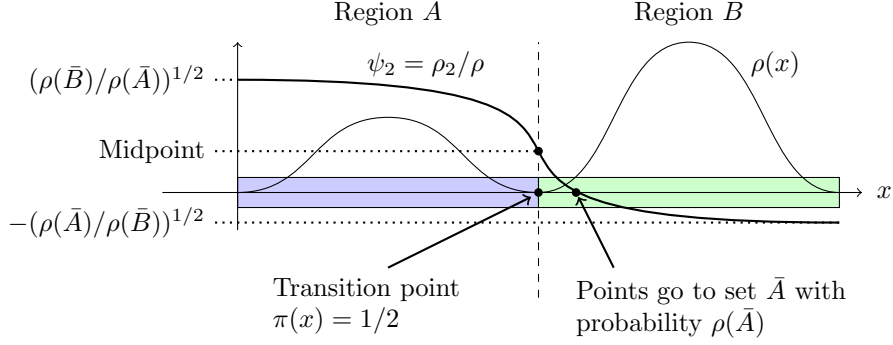

We now discuss in more details the properties of the eigenvectors. If we assume that the relaxation time in basin $\bar{A}$ and $\bar{B}$ is short compared to the reaction rate, i.e., $0<\lambda_2 \ll \lambda_{k}$ for $k>2$, then for $x_0\in\bar{A}$, $x \in \bar{B}$, and $\lambda_k^{-1} \ll t \ll \lambda_2^{-1}$ we have
\[
	\rho (x,t|x_{0},0)\approx\rho(x) + \psi_2(x_0) \rho_2(x) \approx 0 
\]
from which we see that $\rho_2 \propto \rho$ in $\bar{B}$ (and similarly $\bar{A}$), and $\psi_2$ must be approximately constant in $\bar{A}$ (and similarly $\bar{B}$). Moreover, we have $\int \rho_2(x) dx = 0$ from which (with the appropriate normalization):
\begin{equation}
	\rho_2(x) =
	\begin{cases}
		\sqrt{\rho(\bar{B})/\rho(\bar{A})} \; \rho(x), \quad \text{in $\bar{A}$} \\
		- \sqrt{\rho(\bar{A})/\rho(\bar{B})} \; \rho(x), \quad \text{in $\bar{B}$} 		
	\end{cases}	
\end{equation}
Finally we have the general relation
\begin{equation}
	\rho_k(x) = \psi_k(x) \rho(x)
	\label{eq:rhopsi}
\end{equation}
A short proof is provided in the appendix (see Proof 1). The eigenvectors $\rho_2$ and $\psi_2$ are depicted in Fig.~\ref{fig:eig}.

Near the transition region $\psi_2(x) \sim$ constant no longer holds, and in fact, the function $\psi_2$ varies rapidly near the transition region. The committor function $\pi(x)$ is defined as the probability to reach region $B$ before reaching region $A$ starting from $x$. By definition, at the transition region, $\pi(x) = 1/2$. The function $\psi_2$ has a natural interpretation in terms of a committor function (see~\cite{2005JChPh.123m4109R}) and
\begin{equation}
	\pi(x) \sim \frac{\psi_2(x)-\psi_2(a)}{\psi_2(b)-\psi_2(a)} 
	\label{eq:pipsi}
\end{equation}             
See the appendix for a proof of this result (Proof 2).

Therefore, from Eq.~\eqref{eq:pipsi}, the transition region corresponds to $\psi_2(x) = 1/2 \; (\psi_2(A) + \psi_2(B))$, where $\psi_2(A)$ denotes the constant value of $\psi_2$ in $A$ (similarly for $B$). Therefore, at the transition point:
\begin{equation}\label{}
	\psi_2(x) = \frac{1}{2} \; \Big(
	\sqrt{\frac{\rho(\bar{B})}{\rho(\bar{A})} }
	-  \sqrt{\frac{\rho(\bar{A})}{\rho(\bar{B})} }
	\Big)
\end{equation}
The transition point is shown on Fig.~\ref{fig:eig}. The hypersurface corresponding to $\rho_2(x) = 0$ is located nearby but corresponds to a slightly different situation. Points on this hypersurface are not at the transition region but rather are such that they move to basin $\bar{A}$ with probability $\rho(\bar{A})$ and to $\bar{B}$ with probability $\rho(\bar{B})$. (This is true only after a short time $\tau$, with $\tau \gg \lambda_k^{-1}$, $k>2$.) This is shown on Fig.~\ref{fig:eig}.

\bigskip

{\bf Forward and backward rates.} Now that we have derived estimates for $\rho_2(x)$, we can clarify the relation between different rates. To calculate the forward rate, we construct a linear combination $\rho_\text{f}(x,0)$ of $\rho$ and $\rho_2$ such that:
\begin{equation}
	\int \rho_\text{f}(x,0) \; dx = 1, \qquad
	\rho_\text{f}(x,0) \approx 0 \quad \text{for $x \in \bar{B}$.}
\end{equation}
With these conditions we find that the unique solution is:
\begin{equation}
	\rho_\text{f}(x,0) = \rho(x) + \sqrt{ \frac{\rho(\bar{B})}{\rho(\bar{A})} } \rho_2(x)
\end{equation}
Since $\rho$ and $\rho_2$ are eigenvectors:
\begin{equation}
	\frac{\p \rho_\text{f}(x,t)}{\p t}
	= - \lambda_2 \sqrt{ \frac{\rho(\bar{B})}{\rho(\bar{A})} } \; \rho_2(x)e^{-\lambda_{2}t}, \quad
	- \int_{\bar{B}} \frac{\p \rho_\text{f}(x,t)}{\p t} \; dx\bigg|_{t=0}
	= \lambda_2 \, \rho(\bar{B})
\end{equation}
Therefore the forward rate from $A$ to $B$ is equal to $\lambda_2 \, \rho(\bar{B})$. Similarly the rate from $B$ to $A$ is equal to $\lambda_2 \, \rho(\bar{A})$:
\begin{equation}
	\text{rate}_{A \rightarrow B} = \lambda_2 \, \rho(\bar{B}), \quad
	\text{rate}_{B \rightarrow A} = \lambda_2 \, \rho(\bar{A}).
\label{eq19}
\end{equation}

\subsection{Sensitivity of eigenvalues and systematic errors}

\label{sens}

\noindent {\bf Sensitivity of eigenvalue.} We derive a general result regarding the sensitivity of an eigenvalue to perturbations in the matrix entries. We consider a matrix $P$ and assume that we have a full set of eigenvalues. The eigenvalue of interest is $\mu_2$ and:
\begin{equation}
	P - \mu_2 I = R(P)^T (\Lambda(P) - \mu_2 I) \, S(P), \quad
	\text{where $S(P) = [R(P)^T]^{-1}$,}
\end{equation}
and $\Lambda(P)$ is a diagonal matrix containing the eigenvalues. We denote $\bar{P}$ the exact matrix $P$ and consider small deviations $\Delta P = P - \bar{P}$. We assume that in some neighborhood around $\bar{P}$ the matrices $R(P)$, $S(P)$, and $\Lambda(P)$ are differentiable. Calculate the derivative with respect to one of the entries $P_{ij}$:
\begin{equation}
	\frac{\partial P}{\partial P_{ij}}
	= \frac{\partial R(P)^T}{\partial P_{ij}} (\Lambda(P) - \mu_2 I) S(P)
	+ R(P)^T \frac{\partial \Lambda(P)}{\partial P_{ij}} S(P)
	+ R(P)^T (\Lambda(P) - \mu_2 I) \frac{\partial S}{\partial P_{ij}}
	\label{eqdP}
\end{equation}
We denote $e_i$ a row vector such that $[e_i]_j = \delta_{ij}$, $r_2 = R_{2,:}$ (second row of $R$), $s_2 = S_{2,:}$\,. If we evaluate the partial derivative at $P=\bar{P}$, and multiply Eq.~\eqref{eqdP} to the left by $s_2$ and to the right by $r_2^T$ we get:
\begin{equation}
	s_2 e_i^T e_j r_2^T = [s_2]_i \; [r_2]_j =
	\frac{\partial \mu_2}{\partial P_{ij}}
	\label{dmu}
\end{equation}
because $s_2 R^T (\Lambda-\mu_2I) = 0$ and $(\Lambda-\mu_2I) \, S \, r_2^T = 0$. This matrix will be plotted later on, in Fig.~\ref{fig:error7}.

\bigskip

\noindent {\bf Systematic error due to the lag time $\tau$.} We will use this result regarding the sensitivity of $\mu_2$ to study the systematic error or bias using the Markov state model. To simplify the discussion, we will define a fine matrix, which is assumed to accurately capture the continuous dynamics using the Brownian dynamics~\eqref{eq1}:
\begin{equation}
	P^f_{ij} (\tau) = 
	\frac{\int_{x_0 \in V^f_i} \int_{x \in V^f_j} \rho(x,\tau|x_0,0) \rho(x_0) \; dx \; dx_0}{
	\int_{V^f_i} \rho(x_0) \; dx_0
	}
\end{equation}
over some fine states. We consider that $P^f$ gives the exact rate (this can be achieved using states that are fine enough). We could work with $\rho(x,t|x_0,0)$ directly but it is easier to discuss the results using $P^f$. The matrix $P$ can be written in terms of $P^f$:
\begin{equation}
	P_{ij} = \sum_{k \in V_i} \sum_{l \in V_j} \frac{\rho_k}{\sum_{k' \in V_i} \rho_{k'}} \, P^f_{kl} = [\Pi(P^f)]_{kl}
\end{equation}
From $P$, we can reconstruct an approximation $P^c$ of $P^f$ where an entry in $P$ is mapped to a block in $P^c$:
\begin{equation}
	P^c_{kl} = \frac{P_{ij}}{|V_j|} = [{\mathcal I}(P)]_{kl}
\end{equation}
with $k \in V_i$, $l \in V_j$, and where $|V_j|$ is the number of fine states in $V_j$. We have $\Pi \, {\mathcal I} = I$, the identity operator. We will use our result for the sensitivity of $\mu_2$ with $\Delta P = P^c - P^f$.

We now show that every eigenvalue of $P$ is an eigenvalue of $P^c$, which allows us to study the eigenvalues of $P^c$ instead of $P$. Let us define $\tilde{s}_2$ the left eigenvector of $P$ with eigenvalue $\tilde{\mu}_2$. Define $\tilde{s}_2^c$:
\begin{equation}
	[\tilde{s}_2^c]_k = \frac{[\tilde{s}_2]_i}{|V_i|}
\end{equation}
with $k \in V_i$. Then for $l \in V_j$:
\begin{align}
	\sum_k [\tilde{s}_2^c]_k P^c_{kl}
	& = \sum_i \frac{[\tilde{s}_2]_i}{|V_i|} \sum_{k \in V_i} P^c_{kl}
	= \sum_i \frac{[\tilde{s}_2]_i}{|V_i|} \; |V_i| \frac{P_{ij}}{|V_j|} \\
	& = \tilde{\mu}_2 \frac{[\tilde{s}_2]_j}{|V_j|}
	= \tilde{\mu}_2 \, [\tilde{s}_2^c]_l
\end{align}
Therefore $\tilde{s}_2^c$ is a left eigenvector of $P^c$ with eigenvalue $\tilde{\mu}_2$.

Using Eq.~\eqref{dmu}, the sensitivity of the second eigenvalue is:
\begin{equation}
	\frac{\partial \mu_2}{\partial P^f_{ij}} = [s^f_2]_i \; [r^f_2]_j
\end{equation}
where $s^f_2$ is the left eigenvector, and $r^f_2$ is the right eigenvector of $P^f$. We are going to use the following linear approximation:
\begin{equation}
	\tilde{\mu}_2 - \mu_2 \approx
	\sum_{kl} [s^f_2]_k \, [r^f_2]_l \, \big( P^c_{kl} - P^f_{kl} )
\end{equation}   

We start by studying the systematic error at long lag times. Then:
\begin{equation}
	P^f_{kl}
	= \rho(V^f_l) + e^{-\lambda_2 \tau} \; \rho^f_{l,2} \; \langle \psi_2 \rangle^f_k + O(e^{-\lambda_3 \tau})
	\label{Pf}
\end{equation}
In the following we will make the following approximations:\label{assumptions}
\begin{itemize}
	\item $\rho(x)$ is negligibly small near the transition region.
	\item $\psi_2(x)$ is nearly constant away from the transition region.
\end{itemize}
In practice this is not true but the error due to these approximations is typically much smaller than $O(e^{-\lambda_3 \tau})$ and so we will simply ignore it. \footnote{We will not pursue this point further but the analysis suggests that as $\tau \to \infty$ in fact $\tilde{\mu}_2$ does not converge exactly to $\mu_2$ although as explained above this discrepancy is of no practical importance.} As a consequence, from Eq.~\eqref{Pf}, we either have: $P_{kl}^f$ negligible when $l$ is near the transition region or $[r^f_2]_l$ is nearly constant. Denote $[r^f_2]_A$ the value in the left basin and $[r^f_2]_B$ in the right basin. Then:
\begin{equation}
	\sum_l [r^f_2]_l \, P_{kl}^f
	= [r^f_2]_A \sum_{l \in \bar{A}} P_{kl}^f
	+ [r^f_2]_B \sum_{l \in \bar{B}} P_{kl}^f  
    + O(e^{-\lambda_3 \tau})
\end{equation}
In the long lag time assumption, the choice of cells is not important. However, there is one property, which must be satisfied which is that no cell $V_j$ can overlap significantly with both $\bar{A}$ and $\bar{B}$, formally: either $\int_{\bar{A} \cap V_j} \rho(x) dx$ or $\int_{\bar{B} \cap V_j} \rho(x) dx$ must be negligible.
In that case we have the following three possibilities, assuming that $k$ is away from the transition region:
\begin{equation}
	\sum_{l \in V_j} [r^f_2]_l \, P^c_{kl} \approx 
	\begin{cases}
		[r^f_2]_A \sum_{l \in V_j} P_{kl}^f & \text{if $V_j$ falls in $\bar{A}$,} \\
		[r^f_2]_B \sum_{l \in V_j} P_{kl}^f  & \text{if $V_j$ falls in $\bar{B}$,} \\
		0 & \parbox{180pt}{if $V_j$ has support in a region where $\rho$ is \par negligible (transition region).}
	\end{cases}
\end{equation}
with an error of order $O(e^{-\lambda_3 \tau})$. This results from the fact that by construction $P^c_{kl}$ is constant inside $l \in V_j$ and that $\sum_{l \in V_j} P^c_{kl} = \sum_{l \in V_j} P_{kl}^f + O(e^{-\lambda_3 \tau})$. Therefore:
\begin{equation}
	\sum_l [s^f_2]_k \, [r^f_2]_l \, \big( P^c_{kl} - P^f_{kl} )
	= O(e^{-\lambda_3 \tau})
\end{equation} 
At long lag times the error in the eigenvalue is therefore:
\begin{equation}
	\tilde{\mu}_2 - \mu_2 \approx O(e^{-\lambda_3 \tau})
\end{equation}

At short lag times, we have a small systematic error provided the cells are chosen adequately. Let us assume that inside each cell $[r^f_2]_l$ is nearly constant.\footnote{Since $r^f_2$ is an approximation of $\psi_2$ using the fine state discretization (which can be made arbitrarily fine), $r^f_2$ is an approximation of the committor function.} We will therefore denote $[r^f_2]_j$ the value of $[r^f_2]_l$ for $l \in V_j$ in the equation below. The fluctuations of $[r^f_2]_l$ in a given cell are assumed to be of order $\varepsilon$. We also have $[s^f_2]_k = [r^f_2]_k \, \rho(V^f_k)$ [see Eq.~\eqref{eq:rhopsi}], so that:
\begin{align}
	\sum_{kl} [s^f_2]_k \, [r^f_2]_l \, \big( P^c_{kl} - P^f_{kl} )
	& =
	\sum_{ij} [r^f_2]_i [r^f_2]_j \sum_{k \in V_i} \sum_{l \in V_j}
	\rho(V^f_k) \big( P^c_{kl} - P^f_{kl} ) + O(\varepsilon) \\
	& = O(\varepsilon)
\end{align}
from the definition of $P^c$ and $P$.

This suggests the following choice of cell $V_i$ that satisfies our assumption:
\begin{equation}
	V_i = \{ k \, | \, i \varepsilon \le [r^f_2]_k < (i+1) \varepsilon \}
\end{equation}
The interpretation is therefore that the cells provide a fine partitioning based on the iso-surfaces of the eigenvector $\psi_2(x)$ or equivalently the committor function $\pi$, and we can equivalently write:
\begin{equation}
	V_i = \{ k \, | \, i \varepsilon \le \psi_2(x) < (i+1) \varepsilon \}
	\label{eq10}
\end{equation} \label{finecell}
In that case:
\begin{equation}
	\tilde{\mu}_2 - \mu_2 \approx O(\varepsilon)
\end{equation}

This result is consistent with Eq.~(35) in~\cite{Prinz:2011id} which gives a bound on the error that depends on the eigenfunction approximation error when projecting onto the cells. This requirement is very important as constructing cells with small volume is prohibitive in high dimension. However the definition~\eqref{eq10} is much less restrictive in terms of computational cost since it requires refining only along a single direction, given by $\nabla \psi_2$. The dimensionality of the problem has little impact on the number of cells that are required for an accurate calculation. 

The main caveat in this discussion is that computing or even approximating $\psi_2$ or the committor function $\pi$ is very challenging, and consequently this approach remains difficult to use in practice.

\subsection{Statistical errors}

Even though we have outlined methods to reduce the Markovian approximation error, e.g., by refining the spatial discretization, it remains the case in practice that cells need to have a relatively large volume and the committor function is difficult to approximate at best. As a result, a long lag time is more or less necessary for an accurate estimate. We will now discuss how statistical errors vary as a function of the lag time. 

At small lag times, the statistics are typically expected to be reasonably accurate and the statistical error can be made satisfactorily small. At long lag times however the situation worsens. In this section, we will again use Eq.~\eqref{dmu} that expresses the sensitivity of the eigenvalue $\mu_2$ to perturbations in the matrix entries $P_{ij}$, this time due to statistical errors. For this analysis we will assume that the macro states are fine enough, or that $\tau$ is large enough, so that $P(2\tau) = P(\tau)^2$. This simplifies the analysis, but the results mostly carry over to the case where non-Markovian effects are important.

The sensitivity of $\mu_2$ with respect to the matrix entries can be combined with the statistical errors in the entries $P_{ij}$ (\{$P_{i,1}$, \ldots, $P_{i,n_\text{cell}}$\} is a multinomial distribution, \cite{Hinrichs:2007hj,Singhal:2005dk}) to obtain an expression for the statistical error in the form:
\begin{align}
	\sigma^2(\mu_2) & = \frac{1}{n+1} \sum_i \sum_{kl}	  
	\frac{\partial \mu_2}{\partial P_{ik}}
	\frac{\partial \mu_2}{\partial P_{il}}
	\big[ P_{ik} \delta_{kl} - P_{ik} P_{il} \big] \\
	& = \frac{1}{n+1} \big[ \sum_{ik} P_{ik} ([s_2]_i [r_2]_k)^2
	- \sum_i (\sum_k P_{ik} \, [s_2]_i [r_2]_k )^2 \big] \\
	& = \frac{1}{n+1} \sum_i ([s_2]_i)^2 \big[ \sum_k P_{ik} ([r_2]_k)^2
	- (\sum_k P_{ik} \, [r_2]_k )^2 \big] \\
	& = \frac{1}{n+1} \sum_i ([s_2]_i)^2 \; \sigma^2_i(r_2)
\end{align}
where $\sigma_i(r_2)$ is the standard deviation of $r_2$ computed using the probabilities $P_{i,:}$ (row $i$ of $P$).

The rate is given by $\lambda_2 = -\ln(\mu_2)/\tau$ and the relative error can be estimated using $\sigma(\lambda_2) / \lambda_2$:
\begin{equation}
	\frac{\sigma(\lambda_2)}{\lambda_2} \approx 
	\frac{\sigma(\mu_2)}{\mu_2 \tau \lambda_2}	
\end{equation}
We may be interested in minimizing the statistical error, given a computational cost. The cost is proportional to $n\tau/\Delta t$, the number of samples multiplied by the length of the trajectories with lag-time $\tau$. We express the error in the form:
\begin{equation}
	\frac{\sigma(\lambda_2)}{\lambda_2} \approx
	\frac{1}{\sqrt{n \tau / \Delta t} \, \sqrt{\lambda_2 \Delta t}} \frac{\sqrt{\sum_i ([s_2]_i)^2 \sigma^2_i(r_2)}}{e^{-\lambda_2 \tau} \, \sqrt{\lambda_2 \tau}}
	\label{eq:error}
\end{equation}

With our assumption that the states are fine enough, we have that $[s_2]_i \approx \rho_{i,2}$ and $[r_2]_i \approx \langle \psi_2 \rangle_i$. Therefore if we assume that $\lambda_2 \ll \lambda_3$, we have that $r_2$ is nearly constant. Hence $\sigma^2_i(r_2)$ is expected to be small and $\sigma(\lambda_2) / \lambda_2$ can remain bounded even as $\lambda_2 \to 0$.

Again, using the assumption that the states are fine enough, the eigenvectors $s_2$ and $r_2$ of $P$ are independent of $\tau$. At short times, the probability $P_{i,:}$ is concentrated around the diagonal (see Fig.~\ref{fig:error7}). In fact for $\tau = 0$ the matrix is equal to the identity and $\sigma_i(r_2) = 0$. As $\tau$ increases the spread of the entries in $P_{i,:}$ becomes larger and as a consequence $\sigma_i(r_2)$ must increase with $\tau$. The factor $1/(e^{-\lambda_2 \tau}\sqrt{\tau})$ on the denominator in Eq.~\eqref{eq:error} results in an initial increase of the error at small $\tau$, then a plateau is reached, and as $\lambda_2 \tau \gg 1$ the error starts increasing again ($\tau$ is at this point large compared to the relaxation rate of the system). See Fig.~\ref{fig:error7}.

An interpretation is that as $\tau \gg \lambda_3^{-1}$, the system has time to relax within basin $\bar{A}$ or $\bar{B}$. Therefore independent of where the system starts from, the states that were started in basin $\bar{A}$ will be distributed as $\sim \rho(V_i)/\rho(\bar{A})$ while the states in $\bar{B}$ are distributed as $\sim \rho(V_i)/\rho(\bar{B})$. Therefore in that regime, the method degenerates to a direct calculation of the rate where trajectories are initiated in basin $\bar{A}$ (resp.\ $\bar{B}$) and we observe how many transitions to basin $\bar{B}$ (resp.\ $\bar{A}$) occur. For this type of calculation, statistical errors are large whenever $\lambda_2 \tau$ is very small.

\begin{figure}
\centering	
\setlength\figureheight{1.5in} 
\setlength\figurewidth{2in}
    \beginpgfgraphicnamed{tikzfigs/fig1-publisher}%
    \input{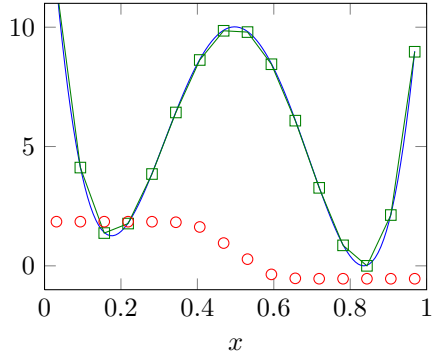}%
    \endpgfgraphicnamed%   
  
\caption{Solid blue line: energy $U(x)$; green line and green squares: fine states used to define the transition matrix $P_{ij}$; red circles: $\psi_2$. The probability density is defined as $\rho = \exp(-U)$.}
\label{fig:error1}
\end{figure}

\begin{figure}
\centering	
\setlength\figureheight{1.25in} 
\setlength\figurewidth{1.6in}
    \beginpgfgraphicnamed{tikzfigs/fig2-publisher}%
    \input{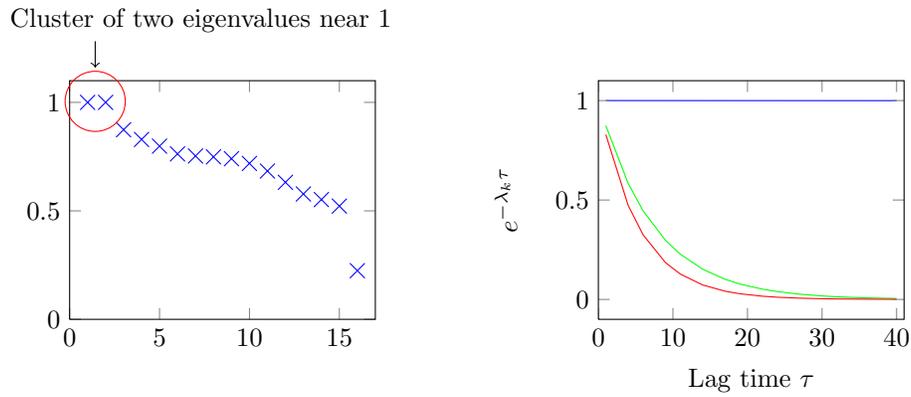}%
    \endpgfgraphicnamed%   

\caption{Eigenvalues of the transition matrix $P$. Left panel: Cluster of eigenvalues near 1. Right panel: Decay of $\exp(- \lambda_k \tau)$ vs $\tau$ for $k=2$, 3, 4.}
\label{fig:error2}            
\end{figure}

\begin{figure}
\centering
    \beginpgfgraphicnamed{tikzfigs/fig6-publisher}%
    \input{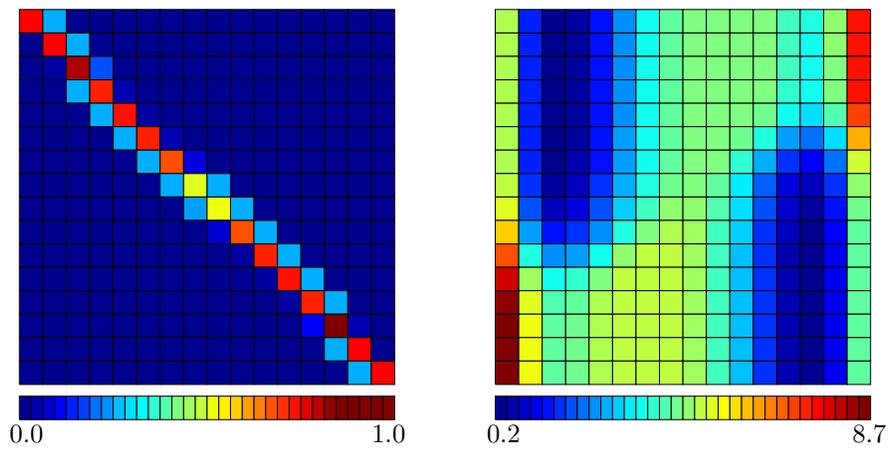}%
    \endpgfgraphicnamed%   
   
\caption{Left panel: matrix $P$; right panel: matrix $-\text{log10}(P^{100})$. The color represents the magnitude of the entries. A color bar is shown at the bottom for reference.}
\label{fig:error7}
\end{figure}

More numerical results will be shown later on, but to illustrate the point above we present a simple example. Fig.~\ref{fig:error1} shows a 1D system with $0 \le x \le 1$. A random walker is moving between discrete states. The probability to attempt a move to the left is 0.25 (same for the right). A Metropolis criterion is used to accept or reject this move (p.~\pageref{metropolis}, and~\cite{frenkel1996understanding} p.\ 27) using:
\begin{equation}
	U(x) = 400 \, \big( 0.98 \, (x-0.2)^4 + (x-0.8)^4 - 1.5 \, (x-0.5)^2 \big)
	\qquad \text{($\beta = 1$)}
\end{equation}
The eigenvalues of the matrix are shown on Fig.~\ref{fig:error2}, along with the decay of $\exp(- \lambda_k \tau)$ as a function of $\tau$ for $k=2$, 3, 4. In Fig.~\ref{fig:error7}, the matrix $P$ is shown along with $-$log10$(P^{100})$. This shows how the matrix $P^{\tau}$ progressively goes from a tri-diagonal form to a rank-2 matrix given by the first two eigenvectors. Recall the assumptions from page~\pageref{assumptions} regarding $\rho(x)$ and $\psi_2(x)$. In Fig.~\ref{fig:error7}, this can be seen from the fact that in the top and bottom portions of the matrix the entries in the columns are nearly constant. In the center of the matrix, entries in the column vary rapidly ($\langle \psi_2\rangle_k^f$ is changing sign) but this is also the region where $\rho$ is small. This can be seen by observing that the columns in the center have small values ($\sim 10^{-4}$). See also Eq.~\eqref{Pf}. The same behavior can be observed in Fig.~\ref{fig:error1} where for $x < 0.4$, $\psi_2$ is constant; for $0.4 \le x \le 0.6$, $\rho$ is small; for $x > 0.6$, $\psi_2$ becomes constant again.

\begin{figure}
\centering	
\setlength\figureheight{1.5in} 
\setlength\figurewidth{2.5in}
    \beginpgfgraphicnamed{tikzfigs/fig8-publisher}%
    \input{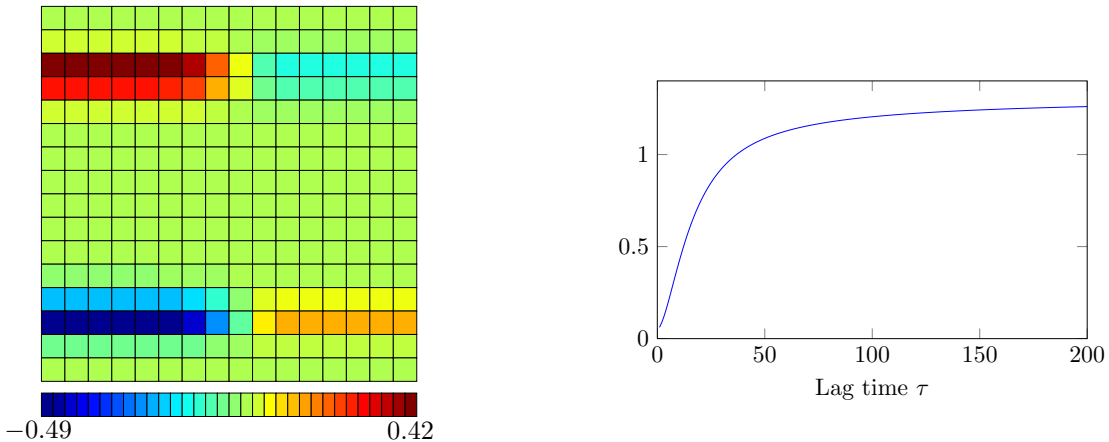}%
    \endpgfgraphicnamed%   
  
\caption{Left panel: sensitivity matrix $\partial \mu_2 / \partial P_{ij} = [s_2]_i [r_2]_j$. Right panel: relative statistical error for the rate estimate as a function of $\tau$ computed using Eq.~\eqref{eq:error}. The statistical error was computed as two standard deviations (95\% confidence interval). The number of sample points used in each state is $10^5/\tau$ [$n \tau / \Delta t = 10^5$ in Eq.~\eqref{eq:error}].}
\label{fig:error8}
\end{figure}

Fig.~\ref{fig:error8} shows the relative statistical error calculated as $2 \sigma(\lambda_2) / \lambda_2$. As $\tau$ increases the statistical error increases because the system has more time to relax in each basin, thereby reducing the computational benefit of using coarse states. 
             
% \begin{figure}
% \centering	
% \setlength\figureheight{1.5in} 
% \setlength\figurewidth{2in}
% \input{tikzfigs/fig4.tikz}
% \caption{Estimated rate using a coarse matrix with only two states. The exact rate is shown with the black line and is equal to 2.86 $10^{-3}$. The entries in the matrix were estimated using $10^4/\tau$ samples in each coarse state. We divided by $\tau$ such that the computational cost is constant as $\tau$ varies. The error bars correspond to 2 standard deviations. 128 runs were calculated to estimate the standard deviation.}
% \label{fig:error4}
% \end{figure}

% \begin{figure}
% \centering	
% \setlength\figureheight{1.5in} 
% \setlength\figurewidth{2in}
% \input{tikzfigs/fig5.tikz}
% \caption{Relative statistical error for the rate estimate as a function of $\tau$. The statistical error was computed as two standard deviations (95\% confidence interval). The number of sample points used in each coarse state is $10^5/\tau$. We used 8 coarse states (see Fig.~\ref{fig:error1}). We used 1024 samples to estimate the standard deviation.}
% \label{fig:error6}
% \end{figure}    

The overall behavior is therefore as follows. At small lag times, we have non-Markovian effects leading to systematic errors. As $\tau$ increases, the terms $e^{-\lambda_k \tau}$, $k>2$, become small when $\lambda_k \tau \gg 1$, such that systematic non-Markovian effects progressively disappear. As this happens, the statistical errors increase because the system has more time to relax in each basin (see Fig.~\ref{fig:error8}). When we reach $\tau \gg \lambda_3^{-1}$ (around $\tau \sim 40$), we see a plateau. Systematic errors are now negligible but the statistical errors are large.

The implication is that it may be difficult in practice to apply these methods accurately. The reasoning above shows that if the cells $V_i$ are not fine enough [see page~\pageref{finecell} and Eq.~\eqref{eq10}], we are caught between systematic errors at small lag times and large statistical errors at large lag times, and a trade-off must be found between these two extremes to maximize the efficiency of the calculation. This issue can be mitigated by choosing cells $V_i$ such that $\psi_2(x) \sim$ constant in each cell, although this can lead to a large computational cost if the number of cells becomes too large or may be intractable if the committor function cannot be well approximated.

The advantage of RTS is that none of these issues are present. Systematic errors are absent and convergence is easy to monitor. There is no lag time $\tau$ that needs to be adjusted to control the accuracy and computational cost. One benefit of MSM though is that it requires only independent sampling inside each macro-state, making the calculation embarrassingly parallel. In contrast, RTS requires a ``global'' convergence of the macro-state weights. This is a necessary trade-off however resulting from the fact that the Markov approximation is not used and therefore transition pathways, in effect, must be computed.             

\section{Some numerical benchmarks for Markov state models and reactive trajectory sampling}

\subsection{Set up of numerical simulations}

To illustrate the performance and numerical behavior of Markov state models (MSM) and the reactive trajectory sampling method (RTS), we present results from some numerical simulations on simple problems in 1D and 2D. They are meant to illustrate the concepts introduced in the previous sections. These results can not really be used as benchmarks of the various methods presented, as conclusions using actual models from bio-chemistry will inevitably be different. Nonetheless they illustrate the key concepts and confirm to some extent the validity of the analysis.

We will consider two types of dynamics. In 1D, we will consider Brownian Dynamics as introduced in Eq.~\eqref{eq1}. Using a time step $\Delta t$, the numerical integrator is given by~\cite{1978JChPh..69.1352E}:
\begin{equation}
x_{n+1} = x_{n} + \nabla D_\text{B} \Delta t - \beta \, D_\text{B} \, \nabla U \, \Delta t
+\sqrt{2\Delta t D_\text{B}} \; W
\end{equation}
where $D_\text{B}$ is the diffusion matrix, $W$ is a vector of independent standard Gaussian random variables (mean $\mu = 0$, standard deviation $\sigma = 1$). A special difficulty of this discretization is that depending on the choice of time step, the walker $x_{n+1}$ may (infrequently) jump over many macro states. In fact there is a finite chance to jump in one step directly from $A$ to $B$. This is an artifact of the discretization. This case is a little artificial but does highlight a limitation of the method, which is that the statistical error typically increases when such ``long-range'' jumps are allowed. In practice however, in the context of molecular systems, this is not an issue as the system typically diffuses slowly and can only cross to neighboring cells during a single time step.

In 2D, the walkers are moving on a 2D Cartesian grid and we consider the Metropolis algorithm (see~\cite{frenkel1996understanding} p.\ 27). \label{metropolis}
The scheme is outlined in algorithm~\ref{algo2}.

\begin{algorithm}
\caption{{\tt Metropolis Algorithm}}   
\label{algo2}
\DontPrintSemicolon
\BlankLine
\tcc{Walkers are moving on a Cartesian grid in dimension $d$.}
\While{more samples are required}{
Move walker to an adjacent position 
$x_{\text{new}}$
with uniform probability of
$1/2^d$ ($d$\nobreakdash-dimensional random walk.)\;
Accept this move with probability 
\[
p=\min\left\{1,\frac{e^{- \beta U (x_{\text{new}})}}
{e^{-\beta U (x_{\text{old}})}}\right\}
\]\;
\tcc{Note that if $U (x_{\text{new}})\le U (x_{\text{old}})$ then $p=1$.}
If move is rejected, stay at $x_{\text{old}}$, otherwise move to $x_{\text{new}}$.\;
}
\end{algorithm}

The exact rate is obtained by computing the eigenvalues of
\[
Q_{ij}=
\mathbf{P} (\text{walker at $j$ the next step} \, | \, \text{walker currently at $i$})
\approx
\rho (x_{j},\Delta t|x_{i},0)
\]
In 2D the system already has discrete states so $Q$ is well-defined. In 1D, we discretize the interval of interest in order to define the fine states $i$.

The coarse transition probability matrix $P$ represents the transition probability between the coarse cells or macro states. The definition is simply as follows:
\begin{equation}
P_{ij} =
\mathbf{P} (\text{walker in $V_{j}$ the next step}\, | \, \text{walker currently in $V_{i}$})
\end{equation}
This matrix can be computed analytically in our examples because everything is low dimensional. For example we can use:
\begin{equation}
	P_{ij} = \sum_{y\in V_{j}} \sum_{x\in V_{i}}
	\mathbf{P}(y|x) \, \mathbf{P}(x|V_{i})
\end{equation}
where $\mathbf{P}(y|x)$ is obtained from $Q$ and $\mathbf{P}(x|V_{i})$ is known exactly. Now given a measurement of $P$, one can compute the second eigenvalue $\mu_{2}$ and let $-\ln (\mu_{2})/\tau $ be an approximation to $\lambda_{2}$.

The RTS algorithm proceeds as follows. In our examples we only have two basins and therefore only use two colors, red and blue. Each time a red particle enters $B$, its color changes to blue, and vice versa. We use Algorithm~{\bf resample} to maintain in each cell or macro state a constant number of walkers of each color. At each step $n$, we calculate
\begin{equation}
	J_n^{\text{b$\to$r}} = \Delta t^{-1} \; \frac{\text{Sum of the weights of all the blue particles who turn red during step $n$}}{\text{Sum of the weights of all the blue particles at step $n$}}
\end{equation}
The rate from $A$ to $B$ is then given by:
\begin{equation}
	\text{rate}_{A \to B} = \lim_{n \to \infty}
	\frac{1}{n} \sum_{i=1}^n J_i^{\text{b$\to$r}}
\end{equation}
The reverse rate, $\text{rate}_{B \to A}$, is obtained similarly.

\subsection{Numerical benchmarks}

{\bf 1D Brownian Dynamics Setup}

The 1D simulation takes place in the domain $\Omega =[-10,10]$. The system parameters are $\beta^{-1} = 0.2$, $D=0.06$ (diffusivity), $\Delta x=0.03$ (fine state discretization), $\Delta t=0.03$, $n_{\text{walkers per cell}} = 10$, and $n_{\text{cell}}=32$ (coarse states). We define: $\bar{A}=[-10,0]$, $A=[-7,-5]$, $\bar{B}=[0,10]$, and $B=[5,7]$. The potential $U(x)$ is
\[
U(x)=\frac{(x+5)^2 (x-5)^{2}}{1000}+3e^{-\frac{x^{2}}{10}}-\frac{x}{10}
\]

\begin{figure}[htbp]
    \center
    \includegraphics[width=0.7\textwidth]{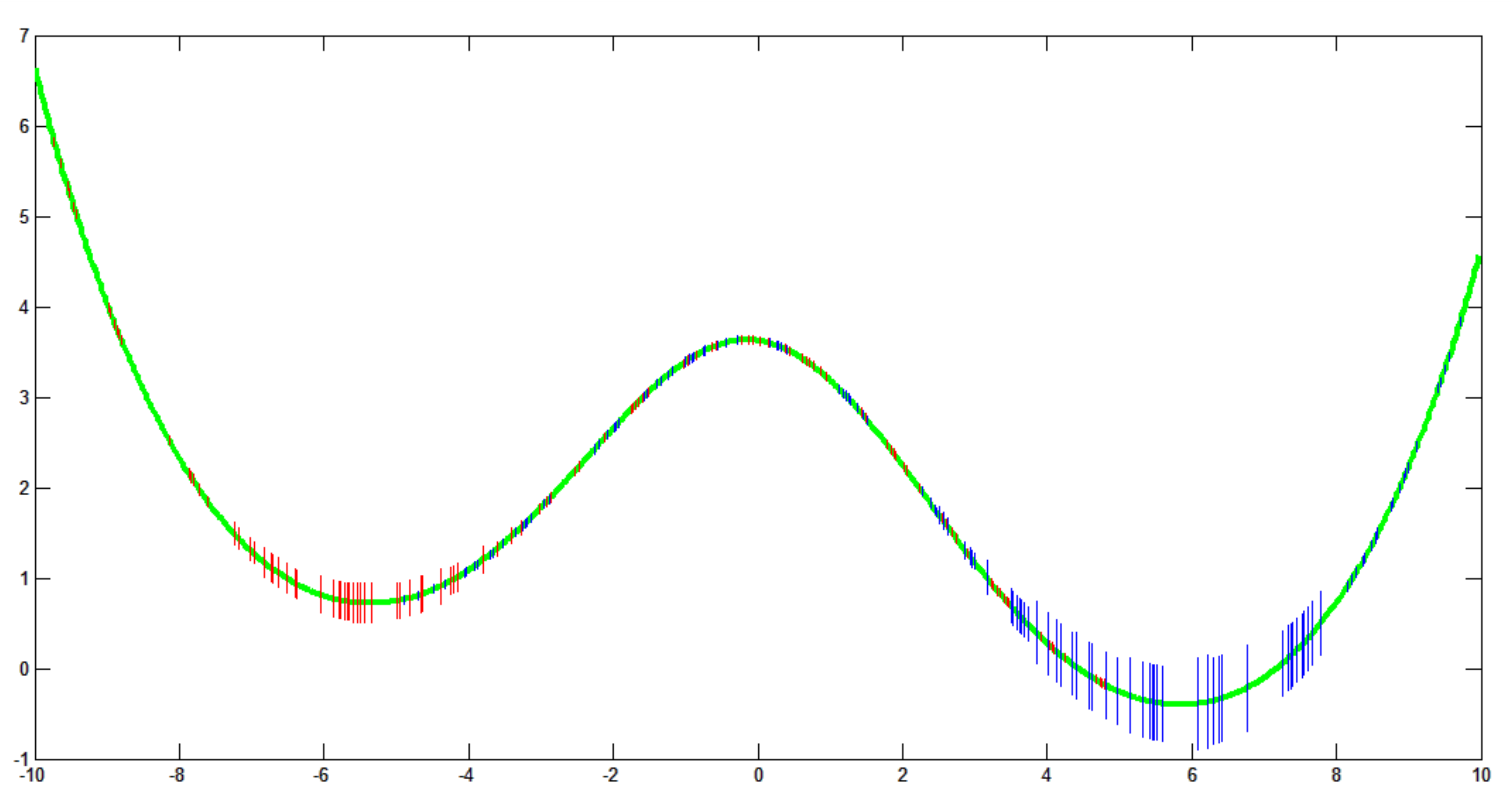}
    \caption{A snapshot of the 1D simulation setup. The green line shows the potential $U(x)$ and the red and blue lines show the positions and weights of the walkers. The weight has large variations and as a result we used a non-linear scaling to determine the length of the short vertical lines with $l=\beta (\log (\alpha w+1)+1)$, where $l$ is the plotted length, $w$ is the weight of the walker, and $\alpha ,\beta$ are constants. (A simple log-scale is not suitable since it would lead to negative vertical line lengths.) }
   \label{fig:1dex}                  %% label for entire figure
\end{figure}

The exact rates computed using the fine states (with spacing $\Delta x$ and time step $\Delta t$) are $\text{rate}_{A\rightarrow B}=1.59\times 10^{-8}$ and $\text{rate}_{B\rightarrow A}=6.70\times 10^{-11}$. The overall setup is illustrated in figure~\ref{fig:1dex}.

The initialization of walkers can impact the convergence time of the rates. A simple choice is to first scatter red walkers uniformly in region $A$ and blue walkers in region $B$ and then assign weight to the walkers proportional to $\rho(x)$. The resulting distribution of red and blue walkers is not the steady-state distribution as the steady-state distribution has red walkers in region $B$ with small weights. However this is a good approximation of this steady-state distribution. The computed rate from the algorithm will not be correct until the distribution of the walkers has converged to steady-state.

\bigskip

\noindent {\bf Discussion}

\nopagebreak

Fig.~\ref{fig1Dc} shows results using the coarse grained matrix and estimating the rate using the second eigenvalue $\mu_2 = e^{-\lambda_2 \tau}$. As predicted by the theoretical analysis of the previous sections we observe two trends. At small lag times $\tau$, the rate is over-predicted because of non-Markovian effects. At longer lag times, the statistical error increases. A good value of the lag time in this example is probably around $\tau \approx 300$, for which the systematic and statistical errors are both reasonably small.

\begin{figure}[htbp]
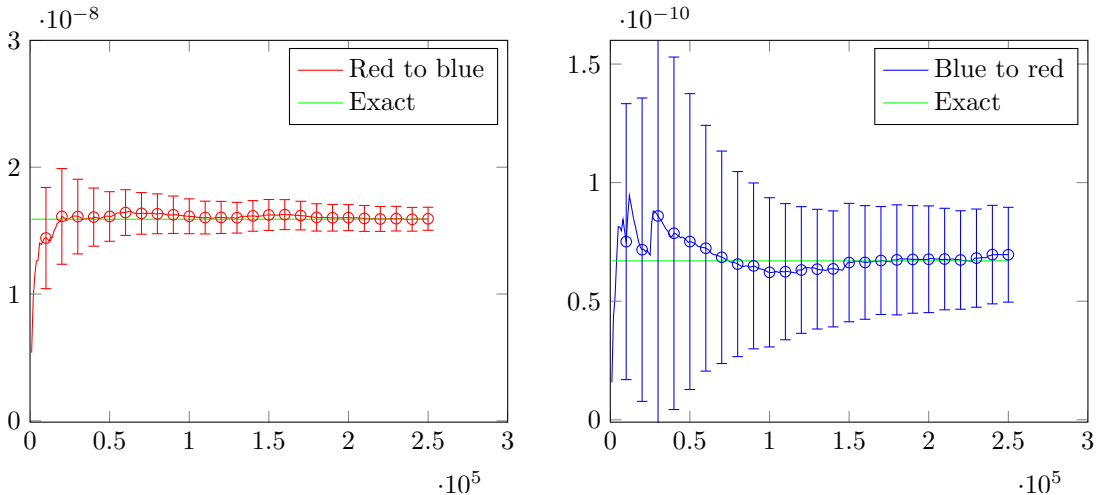

	\center   
	\setlength{\figurewidth}{2.5in}
	\setlength{\figureheight}{2in}
    \beginpgfgraphicnamed{tikzfigs/fig20-publisher}%
    \input{tikzfigs/fig20.tikz}%
    \endpgfgraphicnamed%   

	\quad
    \beginpgfgraphicnamed{tikzfigs/fig21-publisher}%
    \input{tikzfigs/fig21.tikz}%
    \endpgfgraphicnamed%   

	\caption{RTS Rates for the 1D simulation. \label{RTS1D}}
\end{figure}

The rates obtained using the reactive trajectory sampling method (RTS) is shown on Fig.~\ref{RTS1D}. This case is actually relatively difficult as the system is allowed to jump across multiple cells in one time step. As a result there might be some infrequent jump of particles with large weights from basin $A$ to $B$ (or vice versa). This results in larger statistical errors, which leads to a large sampling to reach an acceptable accuracy. The initial bias at small times (underestimation of the rate) is caused by the initial distribution of particles, which needs to equilibrate before correct statistics can be obtained. In practice a more judicious choice of initial distribution for the particles can reduce this initial equilibration.
              
\bigskip

\noindent {\bf 2D Metropolis Method Setup}

\begin{figure}[htbp]
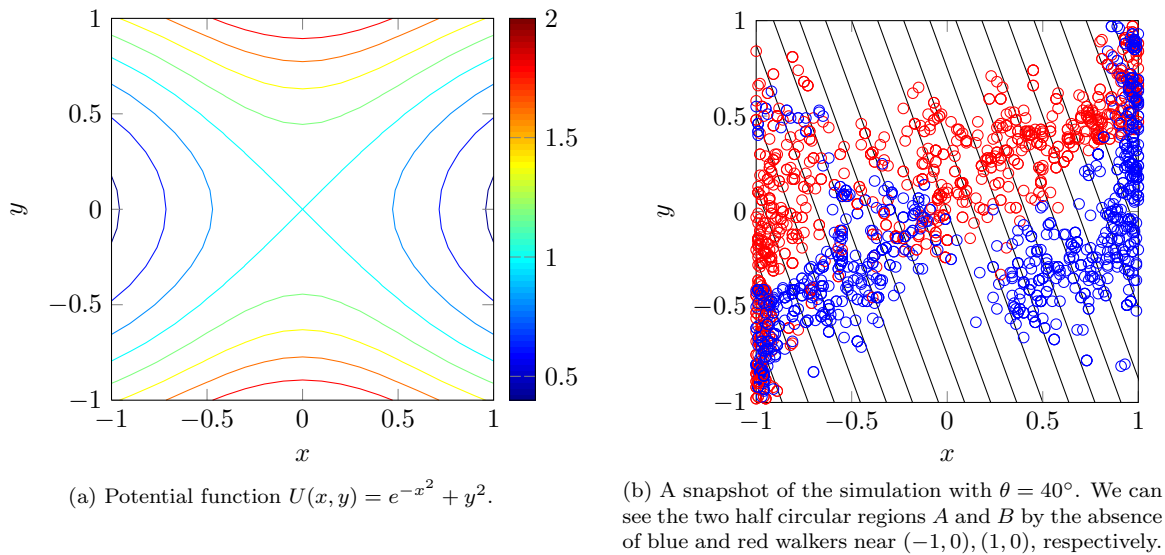

    \center
    \subfloat[Potential function $U(x,y)=e^{-x^{2}}+y^{2}$. \label{subfig:2dex1a} ]{
		\setlength{\figurewidth}{2in}
		\setlength{\figureheight}{2in}
    \beginpgfgraphicnamed{tikzfigs/fig19-publisher}%
    \input{tikzfigs/fig19.tikz}%
    \endpgfgraphicnamed%   

		}
		\quad
    \subfloat[A snapshot of the simulation with $\theta = 40^{\circ}$.
We can see the two half circular regions $A$ and $B$
by the absence of blue and red walkers near $(-1,0),(1,0)$,
respectively. \label{subfig:2dex1b}]{
		\setlength{\figurewidth}{2in}
		\setlength{\figureheight}{2in}
    \beginpgfgraphicnamed{tikzfigs/fig18-publisher}%
    \input{tikzfigs/fig18.tikz}%
    \endpgfgraphicnamed%   

		}
   \caption{Illustration of the 2D simulation setup.}
   \label{fig:2dex1}                  %% label for entire figure
\end{figure}

\begin{figure}[htbp]
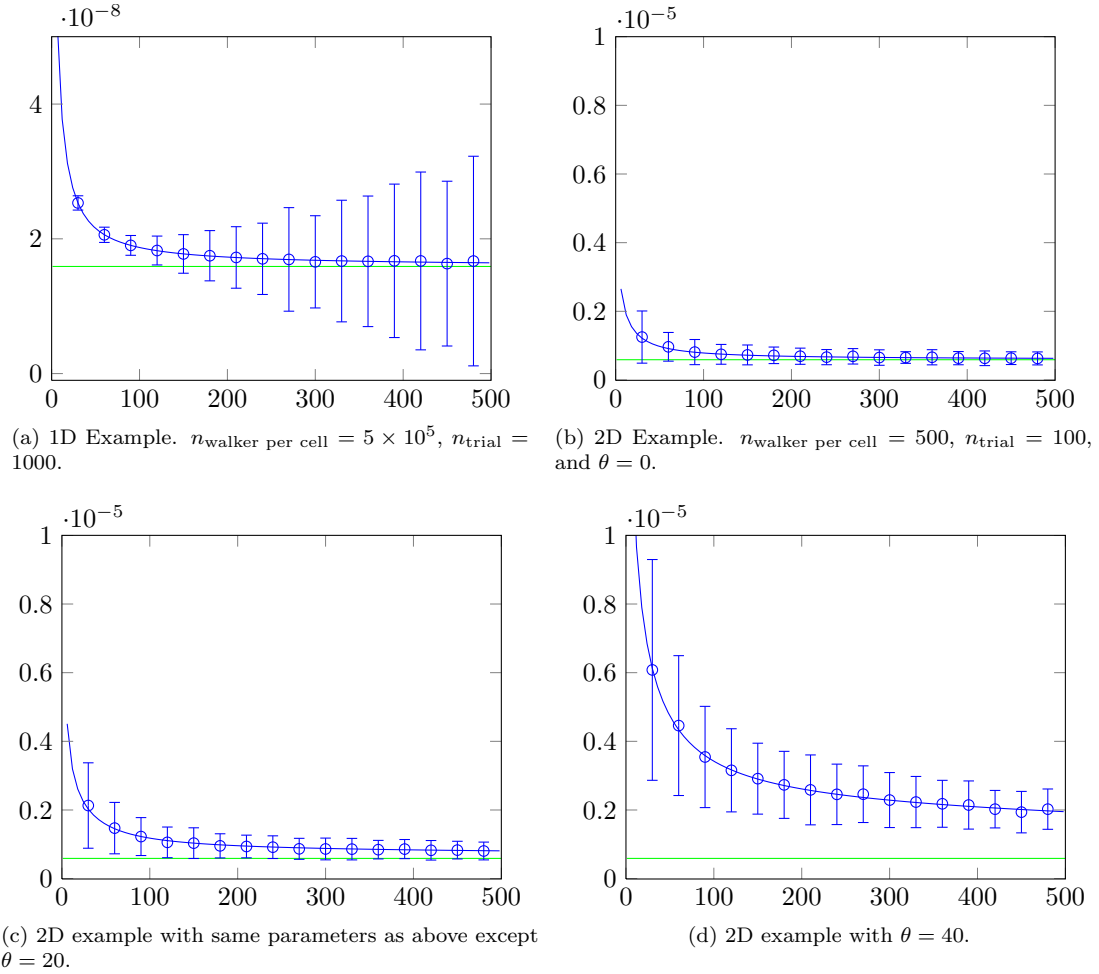

    \center
	\setlength{\figurewidth}{2.3in}
	\setlength{\figureheight}{1.8in}
    \subfloat[1D Example. $n_{\text{walker per cell}}=5\times 10^{5}$, $n_{\text{trial}}=1000$. \label{fig1Dc}
%Data saved in 1Dcmatrixrate.mat and 1DcmatrixRealRate.mat
	]{	
    \beginpgfgraphicnamed{tikzfigs/fig14-publisher}%
    \input{tikzfigs/fig14.tikz}%
    \endpgfgraphicnamed%   

		} \quad
    \subfloat[2D Example.
$n_{\text{walker per cell}}=500$, $n_{\text{trial}}=100$, and
$\theta =0$.  \label{fig2Dc1}
%Data saved in theta0cmatrixrate.mat and theta0MatrixRealRate.mat
	]{   
    \beginpgfgraphicnamed{tikzfigs/fig15-publisher}%
    \input{tikzfigs/fig15.tikz}%
    \endpgfgraphicnamed%   

		}\\
   \subfloat[2D example with same parameters as above
except $\theta =20$.  \label{fig2Dc2}
%Data saved in theta20cmatrixrate.mat and theta20MatrixRealRate.mat
	]{
    \beginpgfgraphicnamed{tikzfigs/fig16-publisher}%
    \input{tikzfigs/fig16.tikz}%
    \endpgfgraphicnamed%   

		}  \quad
   \subfloat[2D example with $\theta =40$.   \label{fig2Dc3}
%Data saved in theta40cmatrixrate.mat and theta40MatrixRealRate.mat
	]{
    \beginpgfgraphicnamed{tikzfigs/fig17-publisher}%
    \input{tikzfigs/fig17.tikz}%
    \endpgfgraphicnamed%   

		}
   \caption{Rates computed with the transition matrix for the coarse grained model. The green line represents the exact rate. The error bars correspond to 2 standard deviations, estimated using $n_{\text{trial}}$ samples. The $x$ axis is the lag time $\tau$.
}
   \label{fig:cmatrate}                  %% label for entire figure
\end{figure}

The 2D simulation takes place in the domain $\Omega=[-1,1]^{2}$. The system parameters are $\beta^{-1} =0.1$, $\Delta x=0.01$, $n_{\text{walkers per cell}} = 100$, and $n_{\text{cell}}=20$. We define the basins as: $\bar{A}=[-1,0]\times [-1,1]$, $A = \{(x,y): \sqrt{(x+1)^{2}+y^2} \le 0.4\}$, $\bar{B}=[0,1]\times [-1,1]$, and $B = \{(x,y): \sqrt{(x+1)^{2}+y^{2}} \le 0.4\}$ (see Fig.~\ref{subfig:2dex1b}). The potential in use is $U(x,y)=e^{-x^{2}}+y^{2}$ and is illustrated in Fig.~\ref{subfig:2dex1a}.

The domain $\Omega$ is divided into $n_{\text{cell}}$ cells where the boundaries are equally spaced lines with a slant angle of $\theta$. Fig.~\ref{subfig:2dex1b} is a snapshot of the simulation with a value of $\theta = 20^{\circ}$ where the boundary of the cells are depicted with the black slanted lines. The time stepping was done using the metropolis algorithm.

If basin $A$ was defined as a vertical slab, the iso-surfaces of the committor function would be vertical lines because of the special form of $U$ in which $x$ and $y$ are basically decoupled. Since $A$ and $B$ are defined as half circles, the choice of $\theta = 0^{\circ}$ makes the cell boundaries close to the iso-surfaces, and thus should be the optimal choice. We will see that as the angle $\theta$ increases the statistical error in the prediction increases.

The exact rates are $\text{rate}_{A\rightarrow B}=\text{rate}_{B\rightarrow A}=5.9\times 10^{-7}$.       
                          
% Directory Dropbox/ernest ryu/data, file paper_plots_walkers.m

\bigskip       

\noindent {\bf Discussion}

\nopagebreak

Results obtained using the coarse grained matrix are shown in Fig.~\ref{fig2Dc1}, \ref{fig2Dc2}, \ref{fig2Dc3}. As before the non-Markovian effects are visible. From the theoretical analysis, we established that non-Markovian effects are reduced when the cells boundary approximate the iso-surface of the committor function and the cells are narrow in the direction $\nabla \pi$. Here the width of the cells is kept fix but we vary the angle of the cells boundaries, which are straight lines. At $\theta = 0$, the lines are relatively good approximation of the iso-surfaces and therefore memory effects are small. As $\theta$ increases, the deviation from the iso-surfaces increases and we observe longer lag times. Note that the cell width is the same so that the degradation does not correspond to a coarsening of the cell but rather a poor choice of their geometry. In complex examples from bio-chemistry we expect that the choice of cells is far from ideal since it is difficult to guess the iso-committor surfaces and therefore the situation in Fig.~\ref{fig2Dc2} and \ref{fig2Dc3} is somewhat representative. We also note that in addition to larger memory effects, the statistical errors increase with $\theta$.

% Directory: Dropbox/Ernest Ryu/data

% Directory containing the data:
% Dropbox/Ernest Ryu/data/theta_x
\begin{figure}[htbp]
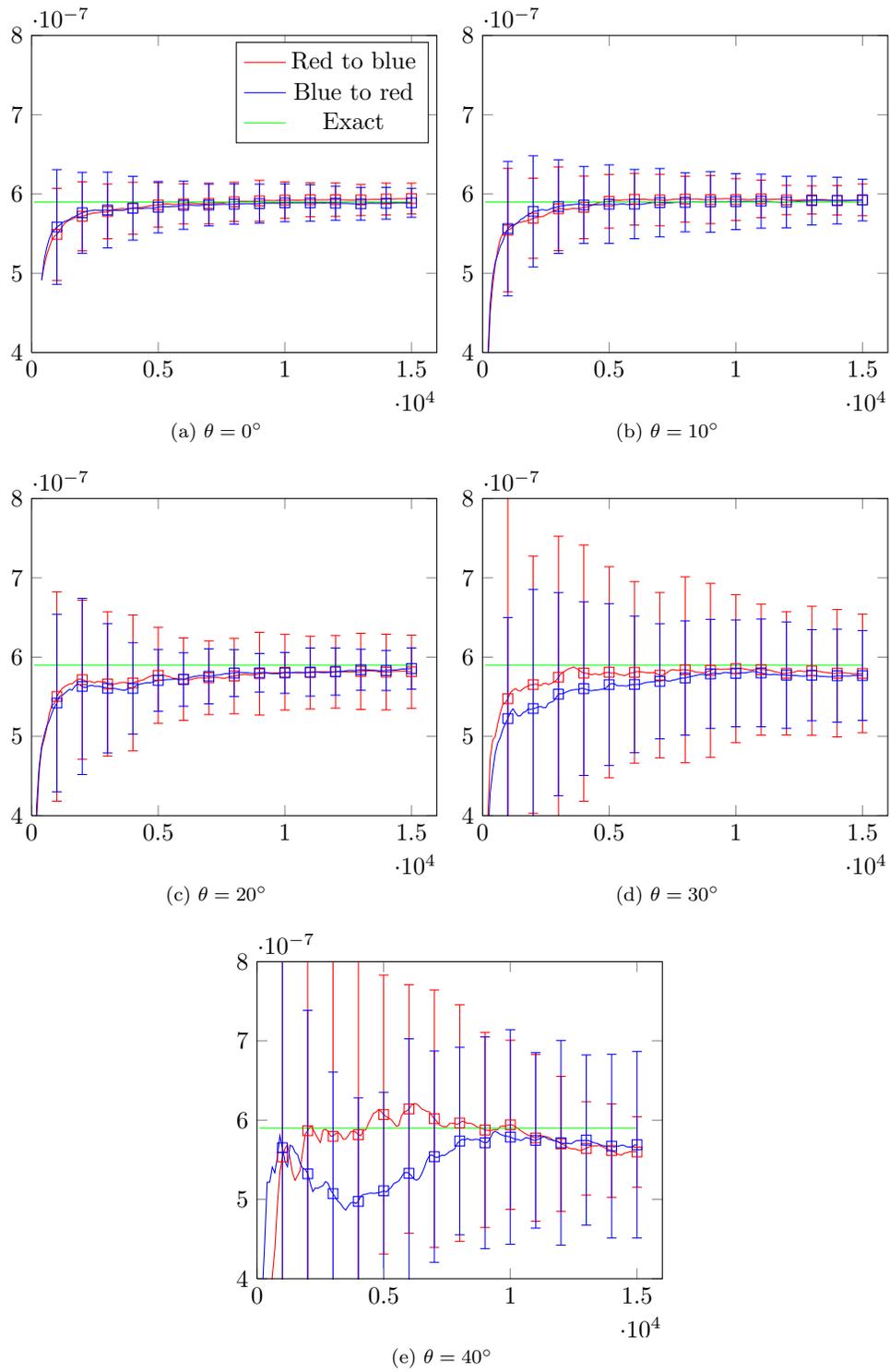

    \center                                                           
	\setlength{\figurewidth}{2.3in}
	\setlength{\figureheight}{1.8in}
    \subfloat[$\theta=0^{\circ}$ \label{fig2Drts1}]{
    \beginpgfgraphicnamed{tikzfigs/fig9-publisher}%
    \input{tikzfigs/fig9.tikz}%
    \endpgfgraphicnamed%   

		}
    \subfloat[$\theta =10^{\circ}$ \label{fig2Drts2}]{ 
    \beginpgfgraphicnamed{tikzfigs/fig10-publisher}%
    \input{tikzfigs/fig10.tikz}%
    \endpgfgraphicnamed%   

		}\\
    \subfloat[$\theta =20^{\circ}$ \label{fig2Drts3}]{
    \beginpgfgraphicnamed{tikzfigs/fig11-publisher}%
    \input{tikzfigs/fig11.tikz}%
    \endpgfgraphicnamed%   

		 }
    \subfloat[$\theta =30^{\circ}$ \label{fig2Drts4}]{
    \beginpgfgraphicnamed{tikzfigs/fig12-publisher}%
    \input{tikzfigs/fig12.tikz}%
    \endpgfgraphicnamed%   

		 }\\
    \subfloat[$\theta =40^{\circ}$ \label{fig2Drts5}]{
    \beginpgfgraphicnamed{tikzfigs/fig13-publisher}%
    \input{tikzfigs/fig13.tikz}%
    \endpgfgraphicnamed%   

		}
   \caption{Rates computed using RTS for different angles $\theta$ for the cells. As $\theta$ increases the method shows no systematic bias but statistical errors increase. The initial relaxation to the exact rate is a result of the initial distribution of particles. The legend is the same for all plots and is only shown for the first plot.} 	
   \label{fig:RTStrate}
\end{figure}

The rates computed using RTS are shown in Fig.~\ref{fig2Drts1}--\ref{fig2Drts5}. We see that as expected the rate always converges towards the correct value. As $\theta$ increases the statistical errors become larger and at $\theta = 40$ the method is no longer very efficient. Even though it is shown, the choice $\theta = 90$ would basically ruin the method since the cells would no longer provide any sampling enhancement. In that case, all the cells would have horizontal boundaries, parallel to $x$.

Fig.~\ref{subfig:2dex1b} illustrates the point made in Section~\ref{opt_cells}. The minimum energy path in this case is the centerline $y=0$. As the cells become tilted particles tend to move away from the centerline. In this case, on the left side, particles tend to move up and on the right side particles tend to move down. Roughly speaking, particles accumulate near the point where $\nabla A$ is orthogonal to the cell boundary. As $\theta$ increases, this point moves away from $y=0$. As a result fewer particles are found near the center lines, even though these are the particles that make the largest contributions to the rate. This translates into infrequent events in which a particle with large weight changes color. At $\theta = 20$ this effect is small but it increases with $\theta$. At $\theta = 90$, we would see very few particles changing color but with a very large weight. The rate would therefore be still correctly calculated, however with a large standard deviation.

\section{Conclusion}

This paper has discussed a number of approaches to calculate reaction rates. We reviewed the reactive flux approach, transition path sampling, transition interface sampling, forward flux sampling, conformation dynamics, Markov state models, non-equilibrium umbrella sampling, and an extension of weighted ensemble Brownian dynamics (renamed reactive trajectory sampling in this paper). These methods differ in their assumptions and computational cost. 

The reactive flux approach is probably the least computationally expensive but requires a fair amount of knowledge about the system (transition pathways, transition region, location of saddle point, etc). Transition interface sampling (TIS) and forward flux sampling (FFS) both rely on a reaction coordinate $\xi$ or an order parameter (variable that increases monotonically with $\xi$). In cases where the mechanism is unknown or multiple pathways are present such methods may be more difficult to apply.

Markov state models (MSM) attempt to alleviate some of these problems by considering a general partitioning of space in macro-states. They can therefore more easily accommodate multiple pathways and complex reaction mechanisms. MSM depends on a lag time $\tau$, equal to the length of trajectories that need to be run. At short times $\tau$, significant non-Markovian effects are present, resulting in a systematic bias, while at long times, statistical errors increase. MSM are ``relatively'' easy to set up and calculate as they involve running a large number of independent short trajectories. It is therefore embarrassingly parallel and can be run efficiently even on loosely connected computers, for example with cloud computing.

Some methods attempt to improve on transition path sampling methods (TIS, FFS) and MSM. They include the non-equilibrium umbrella sampling and weighted ensemble Brownian dynamics, which we called reactive trajectory sampling (RTS) when extended to the case of general macro-states. In this case the assumption of Markovian dynamics can be relaxed. This comes at the cost of a global convergence of statistics across macro-states. For example in RTS and non-equilibrium umbrella sampling, the weights of the macro-states need to be converged, typically using a fixed-point iteration scheme. Although this is typically fast, this does lead in general to a computational cost that is larger than with MSM. In addition, some amount of global communication is required at regular intervals (e.g., for the \verb|resample| algorithm in RTS), although the information that needs to be communicated is minimal. It includes the weights of walkers $w_i$ and macro-state index. Depending on the setup, particle coordinates ($x$) may need to be communicated in some cases.

Techniques exist to construct optimal macro-states or improve an initial guess. They rely for example on computing minimum energy pathways (MFEP). Although there have been many successful implementations, this remains a challenging problem. For example, computing all the MFEPs that make significant contributions to the reaction rate and building macro-states from this data remains a challenge.

Nevertheless, these methods offer promising avenues to calculating reaction rates, transition states and reaction mechanisms. They make excellent use of modern parallel computers as most of these methods involve running a large number of fairly independent trajectories (perhaps with a small amount of global communication required). This puts less pressure on developing software that can run long trajectories using many processors. The efficiency of these methods is independent of the degree of metastability of the system (Perron cluster) and their convergence is often dictated by the shorter mixing time scales inside each metastable basin (at least in an optimal set up of the method).

\medskip

{\bf Acknowledgements.} The authors gratefully acknowledge the work of Jes\'us A.\ Izaguirre and Haoyun ``Michelle'' Feng (University of Notre-Dame) who made several suggestions to improve the manuscript.

\appendix

\section{Technical proofs}

These proofs can be found in classical textbooks such as \cite{Gardiner:1997tb}. We provide them here as they can be helpful to understand some of the results and are also sufficiently simple to be succinctly explained.

\bigskip

\noindent {\bf Proof 1.} We prove that
\begin{equation}
	\rho_k(x) = \psi_k(x) \rho(x)
\end{equation}

Consider simply Eq.~\eqref{eq1} in one dimension (although the proof can be extended to the multi-dimensional case):
\begin{equation}
	dx(t) = A_\text{B}(x) \, dt + \sqrt{B_\text{B}(x)} \; dW(t)
\end{equation}

The function $\rho(x,t|x_0,0)$ satisfies the forward Chapman-Kolmogorov differential equation:
\begin{equation}
	\frac{\p \rho(x,t|x_0,0)}{\p t}     
	= - \frac{d}{dx} (A_\text{B}(x) \rho(x,t|x_0,0))
	+ \frac{1}{2} \frac{d^2}{dx^2} (B_\text{B}(x) \rho(x,t|x_0,0)) 
\end{equation}
From the eigenfunction expansion we therefore have:
\begin{equation}
	-\lambda_k \rho_k(x)     
	= - \frac{d}{dx} (A_\text{B}(x) \rho_k(x))
	+ \frac{1}{2} \frac{d^2}{dx^2} (B_\text{B}(x) \rho_k(x))
	\label{eq6}
\end{equation}   
Similarly the backward equation is satisfied. In this case we write:
\begin{equation}
	\rho(x,t|x_0,0) = \rho(x,0|x_0,-t)
\end{equation}
and
\begin{equation}
	\frac{\p \rho(x,t|x_0,0)}{\p t}
	= -\frac{\p \rho(x,0|x_0,-t)}{\p(-t)}     
	= A_\text{B}(x_0) \frac{d}{dx_0} \rho(x,t|x_0,0)
	+ \frac{1}{2} B_\text{B}(x_0) \frac{d^2}{dx_0^2} \rho(x,t|x_0,0)
\end{equation}      
This leads to:
\begin{equation}
	-\lambda_k \psi_k(x_0)
	= A_\text{B}(x_0) \frac{d}{dx_0} \psi_k(x_0)
	+ \frac{1}{2} B_\text{B}(x_0) \frac{d^2}{dx_0^2} \psi_k(x_0)
	\label{eq7}   	
\end{equation}
For the equilibrium density $\rho$ (which is equal to $\rho_1$):
\begin{equation}
	0 = - \frac{d}{dx} (A_\text{B}(x) \rho(x))
	+ \frac{1}{2} \frac{d^2}{dx^2} (B_\text{B}(x) \rho(x)) \quad \Rightarrow \quad
	A_\text{B} \rho = \frac{1}{2} \frac{d}{dx} (B_\text{B} \rho)
	\label{eq8}
\end{equation}                                         
Consider now $\psi_k(x) \rho(x)$. From Eqns.~\eqref{eq6}, \eqref{eq7}, and \eqref{eq8}, we can prove that:
\begin{equation}
	-\lambda_k \psi_k \rho    
	= - \frac{d}{dx} (A_\text{B} \psi_k \rho_k)
	+ \frac{1}{2} \frac{d^2}{dx^2} (B_\text{B} \psi_k \rho)
\end{equation}
so that:
\begin{equation}
	\rho_k(x) = \psi_k(x) \rho(x)
\end{equation}
$\square$

\bigskip

\noindent {\bf Proof 2.} We prove that
\begin{equation}
	\pi(x) \sim \frac{\psi_2(x)-\psi_2(a)}{\psi_2(b)-\psi_2(a)} 
\end{equation}

We again consider simply a one-dimensional system. Then, the forward Chapman-Kolmogorov equation can be re-written as:
\begin{equation}
	\frac{\p \rho(x,t|x_0,0)}{\p t}     
	= \frac{d}{dx} \Big(-A_\text{B}(x) \rho(x,t|x_0,0)
	+ \frac{1}{2} \frac{d}{dx} (B_\text{B}(x) \rho(x,t|x_0,0)) \Big)
	= -\frac{d}{dx} J(x,t|x_0,0)
\end{equation}
where $J$ is interpreted as a probability flux. Assuming absorbing boundary conditions at $A$ and $B$, the probability to reach $B$ before $A$ is obtained by integrating the flux from 0 to $\infty$:
\begin{equation}
	\pi(x) = \int_0^\infty J(b,t|x,0) \; dt
\end{equation}
assuming that $B=[b,\infty)$. From the definition of the flux:
\begin{equation}
	\pi(x) = \int_0^\infty \Big(
	A_\text{B}(b) \rho(b,t|x,0)
	- \frac{1}{2} \frac{d}{db} (B_\text{B}(b) \rho(b,t|x,0)) 
	\Big) \; dt
\end{equation}   
This time, we use the backward Chapman-Kolmogorov equation in order to obtain derivatives with respect to $x$:
\begin{equation}
	A_\text{B} \, \frac{d\pi}{dx} + \frac{1}{2} \, B_\text{B} \, \frac{d^2\pi}{dx^2}
	= \int_0^\infty \frac{\p J(b,t|x,0)}{\p t} \; dt   
	= -J(b,0|x,0)
\end{equation}
Since for $x \neq b$, $J(b,0|x,0)=0$, we have:
\begin{equation}
	A_\text{B} \frac{d\pi}{dx} + \frac{1}{2} B_\text{B} \frac{d^2\pi}{dx^2} = 0
\end{equation}    
This is the differential equation satisfied by the committor function. This expression can be extended to the multidimensional setting. The boundary conditions are $\pi(a)=0$ and $\pi(b)=1$. The function $\psi_2$ satisfies a slightly different equation:
\begin{equation}
	A_\text{B} \frac{d\psi_2}{dx}
	+ \frac{1}{2} B_\text{B} \frac{d^2 \psi_2}{dx^2} =
	-\lambda_2 \psi_2
\end{equation}
However with $\lambda_2$ very small we can conclude that $\psi_2$ satisfies nearly the same differential equation as $\pi$. From the boundary conditions, we conclude that:
\begin{equation}
	\pi(x) \sim \frac{\psi_2(x)-\psi_2(a)}{\psi_2(b)-\psi_2(a)} 
\end{equation}

$\square$

\bibliographystyle{abbrvnat}

\bibliography{bibliography,library}

\end{document}